\documentclass[12 pt]{article}
\usepackage{amssymb,amsmath,amsfonts,amsthm}
\newcommand\BBP{{\mathbb {P}}}

\newcommand\BBN{{\mathbb {N}}}
\newcommand\BBE{{\mathbb {E}}}

\newtheorem {Lemma}{Lemma}[section]
\newtheorem {Theorem}{Theorem}[section]
\newtheorem {Proposition}{Proposition}[section]

\theoremstyle{definition}
\newtheorem{Definition}{Definition}[section]
\newtheorem{Notation}{Notation}[section]
\newtheorem{Remark}{Remark}[section]

\newcommand\I{{ 1\hspace{-1,2mm}{\mathrm I}}}

\newcommand\no{\noindent}
\newcommand\ssk{\smallskip}

\newcommand\beq{\begin{equation}}
\newcommand\eeq{\end{equation}}

\def\no{\noindent}
\def\ssk{\smallskip}

\def\XZ{(X_i)_{i\in{\mathbb Z}}}
\def\Y{(Y_i)_{i\in{\mathbb N}}}

\def\ZZ{(Z_i)_{i\in{\mathbb Z}}}
\def\EP{(\Omega ,{\mathcal A},{\mathbb P})}
\begin{document}
\title{Deviation inequalities for dependent sequences with applications to strong approximations}
\author{ J. Dedecker, F. Merlev\` ede and E. Rio}
\maketitle

\begin{abstract}
In this paper, we give precise rates of convergence in the strong invariance principle for
stationary sequences of bounded real-valued random variables satisfying weak dependence 
conditions.  One of the main ingredients is a new Fuk-Nagaev type inequality for a class of weakly dependent sequences.  We describe also several classes of  processes to which our results apply.

\medskip

\noindent{\it  MSC2020 subject classifications: }   60F17; 60E15; 37E05. 
\no{\it Running head: }Deviation inequalities and strong approximations.
\end{abstract}

\section{Introduction}

Let $\XZ$ be a strictly
stationary sequence of real-valued random variables (r.v.) defined on a probability space $\EP$, with mean
zero and positive finite variance. Set $S_n = X_1 + X_2 + \cdots + X_n$.  In this
paper, we assume furthermore that the series $\sigma^2=\sum_{k \in {\mathbb Z}} {\rm Cov} (X_0,X_k)$ is convergent (under this assumption $\lim_n n^{-1} {\rm Var} (S_n) = \sigma^2$).  We are interested in obtaining sharp rates in the invariance principle (both in the almost sure sense and in the ${\mathbb L}^2$-sense). Recall that such invariance principles consist in constructing, on a possible larger probability space,  a  sequence $(Z_i)_{i \geq 1}$ of
i.i.d. centered Gaussian variables with variance $\sigma^2$  in such a way  that, setting $T_k = \sum_{i=1}^k Z_i$, 
\beq
\sup_{1 \leq k \leq n} \big | S_k - T_k  \big |  = O (a_n) \text{ a.s. or in ${\mathbb L}^2$}, 
\eeq
where $(a_n)_{n \geq 1}$ is a nondecreasing sequence of positive reals tending to infinity, satisfying $a_n = o ( \sqrt{n} )$. In the independent setting, the almost sure rates are $a_n=o(n^{1/p})$ when ${\mathbb E} (|X_0|^p) < \infty$  for $ p >2$ and are $a_n = O ( \log n)$ when $X_0$ has a finite Laplace  transform in a neighborhood of $0$ (see \cite{KMT,Ma}). Furthermore these rates are optimal according to Breiman \cite{Brei} in the first case and to B\'artfai \cite{Ba} in the second one.  

In the dependent setting, even when the random variables are bounded the almost sure rates can be arbitrarily large. More precisely,  let us consider  the class of irreducible  aperiodic  and positively recurrent Markov chains $(\xi_n)$ with an atom denoted by $A$ (see the definition page 286 in \cite{Bo82}). Let $\tau_A$ be the first return time in $A$, ${\BBP}_A$ be the probability of the chain starting from the atom and ${\BBE}_A$ be the expectation under ${\BBP}_A$. 
  Let then $\pi $ be the unique  invariant distribution, $(\xi_n)$ be the Markov chain starting from  $\pi$, and    $(X_k) $ be the strictly stationary sequence defined by $X_k =f ( \xi_k)$ with $f$ a bounded function.   Theorem 2.2 in \cite{DMR14} asserts that,  for any $p >2$ there exists an irreducible  aperiodic  and positively recurrent Markov chain $(\xi_k)_{k \geq 0}$ with uniform distribution over $[0,1]$ satisfying $\BBP_A ( \tau_A > x) = O ( x^{-p} )  $   and such that for any   absolutely continuous function $f$ on $[0,1]$ with $\pi (f) =0$ and a strictly positive derivative, 
  \[
\limsup_{n \rightarrow \infty} (n \log n)^{-1/p}  \big  | S_n - \sum_{k=1}^n g_k \big | >0 \mbox{ a.s.}
\]
for any stationary and Gaussian centered sequence $(g_k)_{k \in {\mathbb Z}}$ with convergent series of covariances.  This shows that for this type of Markov chains, the rates in the almost sure invariance principle are linked to the moments of the return times in $A$. 

For $p \in ]2,4] $, under the slightly stronger condition $\BBE_A ( \tau_A^p) < \infty$, Cs\'aki and  Cs\"org\H{o}  \cite{CS95} proved  the almost sure rates $a_n= O( n^{1/p }\log n) $ in the strong invariance principle (see their Theorem 2.1).  Our main objective is to extend this result to the case of general stationary sequences of bounded random variables including the case of bounded variation  functions of non irreducible Markov chains.  We shall then consider the case of $\theta$-dependent sequences whose coefficients are defined 
as follows:

\begin{Definition} \label{deftheta} Let 
$ \Gamma_{p,q} = \{  (a_i)_{ 1 \leq i \leq p}   \in {\mathbb N}^p \, : \,  a_1 \geq 1 \text{ and } \sum_{i=1}^p  a_i \leq q  \} $, for $p$ and $q$ positive integers. Let $\XZ$ be a stationary sequence of centered and bounded real-valued random variables and ${\mathcal  F}_0 = \sigma ( X_i, i \leq 0)$. 
 For $k\geq 0$, set
\[
\theta_{X,p,q} (k)  = \sup_{ k_p>k_{p-1}> \ldots > k_2>k_1 \geq k \atop (a_1, \dots, a_p ) \in \Gamma_{p,q}} \Big \Vert  \BBE \Big ( \prod_{i=1}^p X_{k_i}^{a_i} |   {{\mathcal  F}_0 }  \Big  ) - \BBE  \Big ( \prod_{i=1}^p X_{k_i}^{a_i}  \Big )  \Big \Vert_1 \, .
\]
\end{Definition}
These coefficients are suitable for non irreducible Markov chains (see the examples given in Section 4.2). In addition, in case of bounded additive functionals of irreducible  aperiodic  and positively recurrent Markov chains with an atom $A$, for any $p \geq 2$, the condition $\BBE_A ( \tau_A^p) < \infty$ implies that $\sum_{k \geq 1} k^{p-2} \theta_{X,4,4} (k)  < \infty$ (see Section \ref{alphamixing} for more details). Our aim is then to show that for any stationary sequences of bounded random variables satisfying the later weak dependence condition for some $p$ in $ ]2,4]$,   the rate in the almost sure invariance principle is $n^{1/p}$ up to some power of $\log n$ . 

To obtain such rates, a possible approach is to use a martingale approximation and  the Skorokhod embedding theorem. Recall that with this method, the rate cannot be better than $n^{1/4}$ up to some power of $\log n$. The closest results in this direction are given in Wu \cite{Wu} and Doukhan et al \cite{DDM12}. For instance, for $p$ in $]2,4]$, Corollary 3.9 in \cite{DDM12}  provides the rate  $a_n= o( n^{1/p }(\log n)^{1/2+1/p + \varepsilon}) $ for any $\varepsilon>0$ in the almost sure invariance principle    under the condition 
$\sum_{k>0} k^{p-1-2/p} \theta_{X,4,4} ( k) < \infty$, which is  suboptimal.  Still by means of  the Skorokhod embedding theorem, it is  possible to get rates for the  ${\mathbb L}^2$-norm of the error in the invariance principle with suboptimal conditions (see Liu and Wang \cite{LW22}). 

Let $\theta (k) =\theta_{X,4,4} (k) $, for any $k \geq 0$ (recall that the random variables are bounded  in Definition \ref{deftheta}). If $p$ is in $]2,3]$ and  $\sum_{k \geq 1} k^{p-2} \theta (k)  < \infty$, one can obtain the rates $a_n = O ( n^{1/p} )$ (up to some power of $\log n$) both for the almost sure invariance principle and the ${\mathbb L}^2$-norm of the error, by adapting  the proof of \cite[Theorem 2.1]{MR12}.  The main ingredients used in \cite{MR12} are a Fuk-Nagaev type inequality to control the fluctuations and an estimate of the  quadratic cost in the conditional central limit theorem. However, compared to 
inequality \eqref{ineFN1} below, the second term in their  Fuk-Nagaev type inequality cannot be better than  $C n/x^3$. In addition, 
their estimate of the quadratic cost in the conditional central limit theorem cannot be better than $n^{-1/2}$ (this follows from their inequality (A.6)).  These upper bounds induce a limitation of the rates in the  invariance principle  at the level $n^{1/3}$ (up to some power of $\log n$). Recently, the authors proved that for 
$p$  in $]3,4]$ and under the condition $\sum_{k \geq 1} k^{p-2} \theta (k)  < \infty$, the estimate of the quadratic cost in the conditional central limit theorem is of order $n^{-(p-2)/2}$.  This together with the Fuk-Nagaev type inequality stated in Section \ref{SectDI}  are the main ingredients to 
go beyond the rate $n^{1/3}$ in the almost sure invariance principle and to get the rates $n^{1/p}$ (up to some power of $\log n$) for $p \in ]3,4]$. 

Our paper is organized as follows. In Section \ref{SectDI}, we state a Fuk-Nagaev type inequality for partial sums associated with  stationary sequences of bounded  
random variables satisfying the condition $\sum_{k \geq 1} k \theta (k) < \infty$.  In this inequality, the second term is of order $n/x^4$ under the condition 
$\sum_{k \geq 1} k^2 \theta (k) < \infty$.  In Section \ref{SectAP}, we give our main results concerning the rates in the invariance principles (in ${\mathbb L}^2$ and almost surely).  In Section \ref{section3}, we present several classes of examples to which our results apply, including the example of BV observables  of the Liverani-Saussol-Vaienti map (see \cite{LSV}).  An application to rates of convergence in the functional central limit theorem for the quadratic cost associated with the uniform deviation between the Donsker line and the Brownian motion is provided in Section  \ref{SectionRFCLT}. The proof of the Fuk-Nagaev type inequality is given in Section \ref{SectProofFN}. 
 Section \ref{SectASIP} is devoted to the proof of our results concerning the rates in the invariance principles in the non degenerate case, whereas the degenerate case is considered in Section \ref{Sectiondegenerate}. 

In this paper, we shall use the following notations:  $a \ll b$ means that there exists a numerical positive constant $C$ such that  $a \leq C b$, and $X^{(0)}$ means $  X - \BBE (X)$.

\section{Deviation inequalities} \label{SectDI}

\setcounter{equation}{0}

We set $S_0=0$ and $S_n = X_1 + X_2 + \cdots + X_n$ for any positive integer $n$. 
In all the paper, except in Section \ref{Sectiondegenerate},  we  denote  $ \theta_{X,4,4} (k)$ by $ \theta (k) $  
for  all $k$ in $\mathbb N$ and we assume that $\Vert X_0 \Vert_{\infty} =1$. The general case follows by dividing the random variables by $\Vert X_0 \Vert_{\infty}$. 

\begin{Theorem} \label{Fuk-inequality-theta}  Assume that $\sum_{k \geq 1} k \theta (k) < \infty$. 
Let $S_n^* = \max_{0 \leq k \leq n} S_k$.  Then the series 
$\BBE (X_0^2) + 2  \sum_{k \geq 1}\BBE (X_0X_k)$ converges to some nonnegative real $\sigma^2$. 
Moreover,   for any positive real $x$ and any positive integer $n$, 
\beq \label{ineFN1}
\BBP \big ( S_n^* \geq x \big )   \leq c_1 {\bf 1}_{ \{\sigma^2>0\} }  
\Bigl( \frac{   n  \sigma^2 }{x^2} \Bigr)^4  \exp \Big (- \frac{x^2}{  16 n \sigma^2} \Big )  +  c_2 \frac{ n}{ x^{4} }  \Big (  \Theta_1 \Theta_2 + \sum_{k \geq 1} k ( k \wedge x)  \theta (k) \Big )    \, ,
\eeq
where  $\Theta_1 = 1 + \sum_{k\geq 1} \theta (k)$, $\Theta_2 = 1 + \sum_{k\geq 1} k\theta (k)$ 
and $c_1$, $c_2$  are positive numerical constants. 

\end{Theorem}
\begin{Remark}
Recall that $\lim_{n \rightarrow \infty} n^{-1} \BBE (S_n^2) = \sigma^2$ as soon as $\sum_{k \geq 1} \theta(k)< \infty$. 
\end{Remark}

Let $p \in ]3,4[$ and assume that  $\sum_{k \geq 1} k^{p-2} \theta (k) < \infty$. An application of Theorem \ref{Fuk-inequality-theta} gives that for any $\alpha \in ]1/2,1]$ and any $\varepsilon >0$,
$
\sum_{n \geq 1} n^{\alpha p -2}\BBP \big ( S_n^* \geq n^{\alpha} \varepsilon \big )   < \infty $. 
Note that when $p \in ]2,3[$ the same result holds by using Proposition A.2 in \cite{MR12}. 

%\begin{Theorem} \label{Fuk-inequality-thetabis} Let $ \theta (k)  =  \theta_{X,4,4} (k)$. Assume that $\Vert X_0 \Vert_{\infty} \leq M$ and that $\sum_{k \geq 1} k \theta (k) < \infty$. Let $S_n^* = \max_{1 \leq k \leq n} |S_k|$. Then $\sigma^2 = \BBE (X_0^2) + 2  \sum_{k \geq 1}\BBE (X_0X_k)$ and 
%\[\beta_3= \BBE (X_0^3) + 3 \sum_{i \geq 1 }  \{ \BBE ( X_0^2 X_i)  + \BBE ( X_0 X^2_i)  \} + 6 \sum_{u \geq 1} \sum_{v \geq u+1} \BBE ( X_0 X_uX_{v}) \]
%are absolutely convergent. Moreover, for any positive real $x$, 
%\[
%\BBP \big ( S_n^* \geq x \big )   \ll    \exp \Big (- \frac{x^2}{ 4 \kappa_0  n } \Big ) + \frac{ n}{ x^{4} }  \Big ( 1 +  \sum_{k \geq 1} k ( k \wedge x)  \theta_X (k) \Big )  \, ,
%\]
%where $\kappa_0 = 2 \sigma^2 + \frac{ (\beta_3)_+^2}{ 3 \sigma^4}$. 
%\end{Theorem}

\section{Application to strong approximations} \label{SectAP}

\setcounter{equation}{0}

Let $\XZ$ be a stationary sequence of centered and bounded real-valued random variables such that 
$\sum_{k >0}  \theta(k) < \infty$. In this situation, the series $\BBE (X_0^2) + 2  \sum_{k \geq 1}\BBE (X_0X_k)$ is absolutely convergent. In this section, we are interested in strong approximations in the non degenerate case, meaning  that the sum $\sigma^2$ of this series is positive.  We will consider the case $\sigma^2=0$ in Section \ref{Sectiondegenerate}.  In the sequel, we assume that the underlying probability space is rich enough to contain a random variable $\delta$ uniformly distributed over $[0,1]$, independent of the sequence $(X_k)_{k \in {\mathbb Z}}$.  From $\delta$, we construct a sequence $(\delta_i)_{i \in {\mathbb Z}}$ of i.i.d. random variables  uniformly distributed over $[0,1]$ and independent of $(X_k)_{k \in {\mathbb Z}}$.

As explained in \cite{MR12},  having a suitable bound for the quadratic transportation cost in the conditional central limit theorem allows to derive strong approximation results.  Indeed, let us recall the construction given in \cite{MR12} which is inspired from B\'artfai \cite{Ba}.  For $L \in {\mathbb N}$, let $m(L) \in {\mathbb N}$ be such that $m(L)\leq L$. Let 
\begin{equation*} \label{defUkL}
I_{k,L} = ]2^L + (k-1)2^{m(L)} ,  2^L + k 2^{m(L)}] \cap {\mathbb N} \ \text{and}\ U_{k,L} = \sum_{i\in I_{k,L}} X_i \, ,
\, k \in \{1, \cdots, 2^{L-m(L)} \} \, .
\end{equation*}
For $k \in \{1, \cdots, 2^{L-m(L)} \} $, let $V_{k,L}$ be the
$N (0, \sigma^2 2^{m(L)} )$-distributed random variable
defined from $U_{k,L}$ via the conditional quantile transformation,
that is \beq \label{defVkN} 
V_{k,L} = \sigma 2^{m(L)/2} \Phi^{-1}
(\widetilde F_{k,L}( U_{k,L}- 0 ) + \delta_{2^L + k 2^{m(L)}} (\widetilde
F_{k,L}( U_{k,L} ) - \widetilde F_{k,L} (U_{k,L} - 0) ) )  \, ,
\eeq where
$\widetilde F_{k,L}:= F_{U_{k,L} | {\mathcal F}_{2^{L}+ (k-1)2^{m(L)}}}$ is the d.f. of 
$P_{U_{k,L} | {\mathcal F}_{2^L + (k-1)2^{m(L)}
}}$ (the conditional law of $U_{k,L}$ given $\mathcal{F}_{2^L +
(k-1)2^{m(L)} }$)     and $\Phi^{-1}$ is the inverse of the
standard Gaussian distribution function $\Phi$. Since $\delta_{2^L + k 2^{m(L)}}$ is
independent of  ${\mathcal F}_{2^L + (k-1)2^{m(L)}}$, the random variable $V_{k,L}$
is independent of ${\mathcal F}_{2^L + (k-1)2^{m(L)}}$, and has the Gaussian distribution
$N (0, \sigma^2 2^{m(L)})$ (see \cite[Lemma F.1]{Rio17}. By induction on $k$, the random variables $(V_{k,L})_k$
are mutually independent and independent of ${\mathcal F}_{2^L}$. 
In addition
 \begin{eqnarray} \label{L2wasser}\BBE
(U_{k,L} - V_{k,L})^2 & = & \BBE \int_0^1 \big ( F^{-1}_{U_{k,L} |
{\mathcal F}_{2^L + (k-1)2^{m(L)} }}(u)- \sigma 2^{m(L)/2}\Phi^{-1}(u)
\big )^2 du \nonumber \\ & : = &  \BBE  \big ( W_2^2 (P_{U_{k,L} | {\mathcal
F}_{2^L + (k-1)2^{m(L)} }} , G_{\sigma^2 2^{m(L)} }) \big ) \, ,
\end{eqnarray} 
where $G_{\sigma^2 2^{m(L)} }$ is the Gaussian distribution
$N (0, \sigma^2 2^{m(L)})$. 

This construction together with Theorem \ref{Fuk-inequality-theta}  and estimates of the quadratic conditional cost in the central limit theorem given in \cite{DMR22} are  the main ingredients to get   the next strong approximations results.

\begin{Theorem} \label{thSAPfaible}  $ $
\begin{enumerate}
\item[(i)]   Assume that $ \theta(k) = O ( k^{1-p} )$  if $p \in ]2,3[ \cup ]3,4[$ and 
$\sum_{k >0} k^{p-2} \theta(k) < \infty$ if $p \in \{3,  4\}$.  Then, one can construct a sequence of i.i.d.  Gaussian random variables $(Z_i)_{i \geq 1}$ centered and with variance $\sigma^2$, such that, setting $T_k = \sum_{i=1}^k Z_i$, we have
\[
  \Big \Vert \sup_{k \leq n} |S_k-T_k|   \Big \Vert_2  =  O (n^{1/p} (\log n)^{ (1/2 - 1/p) }  ) \, .
\]
\item[(ii)]   Let $p \in \{3,  4\}$ and assume that $ \theta (k) = O ( k^{1-p} )$. Then, one can construct a sequence of i.i.d.  Gaussian random variables $(Z_i)_{i \geq 1}$ centered and with variance $\sigma^2$, such that, setting $T_k = \sum_{i=1}^k Z_i$,
\[
  \Big \Vert \sup_{k \leq n} |S_k-T_k|   \Big \Vert_2   =  O (n^{1/p} (\log n)^{ 1/2  } )  \, .
\]
%\item[(iii)]   Let $p \in \{3,  4\}$ and assume that $ \sum_{k >0} k^{p-2} \theta (k) < \infty$.  Then, one can construct a sequence of i.i.d.  Gaussian random variables $(Z_i)_{i \geq 1}$ centered and with variance $\sigma^2$, such that, setting $T_k = \sum_{i=1}^k Z_i$,
%\[
%  \Big \Vert \sup_{k \leq n} |S_k-T_k|   \Big \Vert_2  = O (n^{1/p} (\log n)^{ (1/2 - 1/p) }  )   \, .
%\]
\end{enumerate}
\end{Theorem}

We now give almost sure rates for strong approximations. 

\begin{Theorem} \label{thSAP}  $ $ 
\begin{enumerate}
\item[(i)] Let $p \in ]2,3[ \cup ]3,4[$ and assume that $ \sum_{k >0} k^{p-2} \theta (k) < \infty$.   Then,  one can construct a sequence of i.i.d.  Gaussian random variables $(Z_i)_{i \geq 1}$ centered and with variance $\sigma^2$, such that, setting $T_k = \sum_{i=1}^k Z_i$,
\[
\sup_{k \leq n} |S_k-T_k| = o (n^{1/p} (\log n)^{1/2-1/p}) \text{ a.s.}
\]
\item[(ii)] Assume that $ \theta (k) = O ( k^{1-p} )$  if $p \in ]2,3[ \cup ]3,4[$ and 
$\sum_{k >0} k^{p-2} \theta (k) < \infty$ if $p \in \{3,  4\}$.   Then,  for any $\eta >1/2$, one can construct a sequence of i.i.d.  Gaussian random variables $(Z_i)_{i \geq 1}$ centered and with variance $\sigma^2$, such that, setting $T_k = \sum_{i=1}^k Z_i$,
\[
\sup_{k \leq n} |S_k-T_k| = o (n^{1/p} (\log n)^{\eta}) \text{ a.s.}
\]
\item[(iii)]    Let $p \in \{3,  4\}$ and assume that $ \theta (k) = O ( k^{1-p} )$.  Then,  for any $\eta >1/2+1/p$, one can construct a sequence of i.i.d.  Gaussian random variables $(Z_i)_{i \geq 1}$ centered and with variance $\sigma^2$, such that, setting $T_k = \sum_{i=1}^k Z_i$,
\[
\sup_{k \leq n} |S_k-T_k| = o (n^{1/p} (\log n)^{\eta}) \text{ a.s.}
\]
\end{enumerate}
\end{Theorem}

\begin{Remark} Let $(\beta (k))_{k \geq 1} $ be  the  sequence of the usual $\beta$-mixing coefficients. 
According to the definition of $\beta(k)$ given page 147 in \cite{MPU},  $\theta( k) \leq 2 \Vert X_0\Vert_\infty^4\, \beta(k)$. Hence, by Item (b) of Theorem 2.2 in \cite{DMR14},  for any $p >2$ there exists a stationary Markov chain $(X_k)_{k \geq 0}$ with uniform distribution such that $\theta( k) \leq C k^{1-p}$ and  
\[
\limsup_{n \rightarrow \infty} (n \log n)^{-1/p}  \big  | S_n - \sum_{k=1}^n g_k \big | >0 \mbox{ a.s.}
\]
for any stationary and Gaussian centered sequence $(g_k)_{k \in {\mathbb Z}}$ with convergent series of covariances.  Consequently,  the rates in Theorem \ref{thSAP}  (and then also in Theorem \ref{thSAPfaible}) are optimal up to 
a power of $\log n$. 
\end{Remark}

\section{Examples and applications} \label{section3}

\setcounter{equation}{0}

\subsection{$\alpha$-mixing sequences}\label{alphamixing}
Let $(\Omega, {\mathcal A}, \BBP)$ be a probability space and let $ {\mathcal U}$ and $ {\mathcal V}$
 be two $\sigma$-algebras of $ {\mathcal A}$.  The strong mixing coefficient $ \alpha (\mathcal{U},\mathcal{V})  $ between these $\sigma$-algebras is defined as follows: 
 \[
  \alpha (\mathcal{U},\mathcal{V})  =\sup \{ | {\ \mathbb{P}} (U\cap
V) - {\mathbb{P}} (U) { \mathbb{P}} (V)|: U \in \mathcal{U}, V\in 
\mathcal{V} \} \, .
 \]
Next, for a stationary sequence $(Y_i)_{i \in {\mathbb{Z}}}$ of random variables with values in a Polish space $S$, define its 
 strong mixing (or $\alpha$-mixing) coefficients of order $4$  as follows:  Let  
\[
 \alpha_{\infty, 4}(n) =  \sup_{i_4 > i_3 >i_2 > i_1 \geq n }  \alpha ({\mathcal F}_{0},  \sigma (Y_{i_1}, Y_{i_2}  , Y_{i_3}, Y_{i_4}) )   \, ,
 \] 
 where ${\mathcal F}_0= \sigma(Y_i, i \leq 0 )$.  As page 146 in \cite{MPU}, these coefficients can be rewritten in the following form:  Let $B_1$ be the class of  measurable functions   from $S^4$ to ${\mathbb R} $ and bounded by one. Then
 \[
  \alpha_{\infty, 4}(n)  = \frac{1}{4}   \sup_{f \in B_1}  \sup_{i_4 > i_3 >i_2 > i_1 \geq n }   \big   \Vert  \BBE ( f( Y_{i_1}, Y_{i_2}  , Y_{i_3}, Y_{i_4})  | {\mathcal F}_0 ) -  \BBE ( f( Y_{i_1}, Y_{i_2}  , Y_{i_3}, Y_{i_4}) )  \big \Vert_1 \, .
 \]
 Let now $f $ be a bounded measurable numerical  function and let $X_k= f(Y_k) - \BBE(f(Y_k))$. Then Theorems \ref{thSAPfaible} and \ref{thSAP} apply to the partial sum   $S_n = \sum_{k=1}^n X_k$, replacing the conditions on   $\theta (k)$ by the same conditions on $\alpha_{\infty, 4}(k)$.  
%For instance, the following corollary holds.
 %\begin{Corollary} \label{coralphamixing}
% Let $(Y_k)_{k \in {\mathbb Z}}$ be a stationary sequence of random variables  with values in a Polish space and such that  $\alpha_{\infty, 4} (k)=O(k^{-3})$.
 
\medskip

Consequently, our  results apply to irreducible stationary $S$-valued Markov chains $(Y_i)_{i  \geq 0}$ with  invariant probability $\pi$ and transition kernel $P(\cdot, \cdot)$ satisfying  the following minorization condition: there exists a positive integer $m$ such that 
\[
P^m ( x, A ) \geq s(x) \nu ( A) \quad \mbox{for $x\in S$, $A \in {\mathcal B}(S)$}\, ,
\]
where $s$ is a measurable function with values in $[0,1]$ such that $\pi(s) >0$  and $\nu$ is a probability measure absolutely continuous with respect to $\pi$ (see \cite[Chapter 2]{Nu84} and \cite[Th. 9.2.15]{DMPS18} for the fact that the invariant probability
measure is a maximal irreducibility measure).

Indeed, let us explain  how the strong mixing coefficients of  the chain $(Y_i)_{i  \geq 0}$ can be computed (we simply denote by
$\alpha(k)$ these coefficients, because in this Markovian setting,
$\alpha(\sigma(Y_0), \sigma(Y_k))= \alpha(\sigma(Y_0), \mathcal{G}_k)$, where $\mathcal{G}_k= \sigma(Y_i, i \geq k)$).
Let $Z_k=Y_{km}$. It follows that  $(Z_k)_{k \geq 0}$ is an  irreducible stationary Markov chain satisfying the  minorization condition with $m =1$ and then the conditions of \cite[Proposition 9.7]{Rio17}. Let $\alpha_Z(k)$ denote the $\alpha$-mixing coefficients of $(Z_i)_{i  \geq 0}$. According to \cite[Page 165]{Rio17} (see also \cite{Bo82}),  if $\tau$ is one of the return times (i.e. the difference between two regeneration times) of the extended chain constructed from $(Z_i)_{i  \geq 0}$, then for $q\geq 1$, $\sum_{k \geq 1} k^{q-1} \alpha_Z(k) < \infty$ iff ${\mathbb E}(\tau^{q+1}) <\infty$. Now, since 
$
\alpha (k) \leq \alpha_Z ([k/m])$,
we infer that ${\mathbb E}(\tau^{q+1}) < \infty$ also implies that $\sum_{k \geq 1} k^{q-1} \alpha(k) < \infty$. In particular, if $\tau$ is such that ${\mathbb E}(\tau^{p}) < \infty$, then $\sum_{k \geq 1} k^{p-2} \alpha(k) < \infty$,  and Theorems \ref{thSAPfaible} and \ref{thSAP} apply to $S_n=\sum_{k=1}^n ( f(Y_k)-{\mathbb E}(f(Y_k)))$.

In their paper, Merlev\`ede and Rio \cite{MR15} proved the almost sure invariance principle with rate $O(\log n)$ when the Markov chain is geometrically ergodic and the minorization condition holds with $m=1$. This last condition allows the use of the 
regeneration technique. However in some situations one can only prove that the  minorization condition  holds for $m >1$ (see for instance Exemple 2.3 (f) in Nummelin \cite{Nu84}), or even that it cannot hold with  $m=1$ 
(even with the additional assumption that there exists a joint density for $(X_0,X_1)$, see Kendall and Montana \cite{KM}). 
Nevertheless, as explained before,  for this class of Markov chain we do not need to assume $m =1$ since our conditions are only expressed in terms of the $\alpha$-mixing coefficients of the chain. 
%More precisely, defining $(Y_k)_{k \geq 0}$ by $Y_k = X_{km}$, one can see that $(Y_k)_{k \geq 0}$ is an  irreducible %stationary Markov chain with the minorization condition which holds with $m =1$  
%\[
%\beta_X (k) \leq \beta_Y ([k/m])  \, .
%\]

\subsection{$\alpha$-dependent sequences}\label{taualpha} 

We start by  recalling the definition of the $\alpha$-dependence coefficients as considered in  \cite{DGM10}. 
\begin{Definition}\label{defalpha}
For any random variable $Y=(Y_1, \cdots, Y_k)$ with values in
${\mathbb R}^k$ and any $\sigma$-algebra ${\mathcal F}$, let
$\alpha({\mathcal F}, Y)= \sup_{(x_1, \ldots , x_k) \in {\mathbb R}^k}
\big \| \BBE \bigl( \, \prod_{j=1}^k (\I_{Y_j \leq x_j})^{(0)} \mid {\mathcal F} \bigr)^{(0)} \big\|_1$,
where we recall that $Z^{(0)}$ means $Z-{\mathbb E}(Z)$. 
For the sequence ${\bf Y}=(Y_i)_{i \in {\mathbb Z}}$, let ${\mathcal F}_0 = \sigma ( Y_i, i \leq 0)$,
\begin{equation}
\label{defalpha} \alpha_{{\bf Y},4}(0) =1 \text{ and }\alpha_{{\bf Y},4}(n) = \max_{1 \leq l \leq
4} \ \sup_{ n\leq i_1\leq \cdots \leq i_l} \alpha({\mathcal F}_0,
(Y_{i_1}, \dots, Y_{i_l})) \text{ for $n>0$}\, .
\end{equation}
\end{Definition}
Let BV$_1$ be the space of bounded variation functions $f$ such that $\|df\|\leq 1$, where $\|df\|$ is the variation norm on ${\mathbb R}$ of the measure $df$. As mentioned in \cite{DGM10}, $\alpha({\mathcal F}, Y)$ can also be defined by
$$
\alpha({\mathcal F}, Y)= \sup_{(f_1, \ldots , f_k) \in \text{BV}_1}
\Big \| \BBE \Big(\prod_{j=1}^k f_j(Y_j )^{(0)} \Big | {\mathcal F} \Big)^{(0)} \Big\|_1 \, ,
$$
It follows that, if $f$ is a bounded variation function such that $\|df\|\leq C$, and $X_k=f(Y_k)-{\mathbb E}(Y_k)$, then 
$\theta (k) \leq C^4 \alpha_{{\bf Y},4}(k)$. 
Then Theorem \ref{thSAPfaible} and \ref{thSAP} apply to the partial sum   $S_n = \sum_{k=1}^n X_k$, replacing the condition on   $\theta (k)$ by the same conditions on $\alpha_{{\bf Y},4}(k)$.  

%\begin{Corollary} \label{coralpha}
%Let $f $ be a bounded variation (BV)  function and $X_k= f(Y_k) - \BBE(f(Y_k))$ where $(Y_k)_{k \in {\mathbb Z}}$ is a %stationary sequence of real-valued random variables. Let $S_n = \sum_{k=1}^n X_k$. If
%$\sum_{k \geq 1} k^2 \alpha_{4, {\bf Y}} (k) < \infty$,
%then  $  W_2 (P_{S_n / {\sqrt n} } ,  G_{\sigma^2}) \ll  n^{-1/2}$.  
%\end{Corollary}
From this result and proceeding as in \cite[Section 3]{MR12}, we can derive rates in the strong approximation results for the partial sums associated with  BV observables of the LSV map. More precisely, for  $\gamma \in ]0,1[$, let  $T_\gamma$ defined from $[0,1]$ to $[0,1]$ by
\begin{equation*}
T_\gamma(x)=
\begin{cases}
x(1+ 2^{\gamma}x^{\gamma}) \quad \text{ if $x \in [0, 1/2[$} \\
2x-1 \quad \quad \quad \ \ \text{if $x \in [1/2, 1]$} \, .%
\end{cases}
\end{equation*}
This is the so-called LSV  \cite{LSV} map with parameter $\gamma$. Recall, that there exists a unique $T_\gamma$-invariant measure $\nu_\gamma$ on $[0, 1]$, which is absolutely
continuous with respect to the Lebesgue measure with positive density denoted by $h_\gamma$.  From \cite[Prop. 1.17]{DGM10}, we know that the coefficients $\alpha_{{\bf Y},4}(k)$ of the  Markov chain associated with $T_\gamma$ are exactly of order $1/k^{(1-\gamma)/\gamma}$. Consequently, if  $f$ is a BV observable, we get that:
\begin{itemize}
\item For any $\gamma \in ] 0, 1/2[$ and any $\varepsilon >0 $, one can construct on  the probability space $([0,1]\times[0,1], \nu_{\gamma}\otimes \lambda)$ a sequence of i.i.d.  Gaussian random variables $(Z_i)_{i \geq 1}$ centered and with variance $\sigma^2$, such that
\begin{equation}\label{AS}
\sup_{k \leq n} \left |\sum_{i=1}^k (f\circ T_{\gamma}^{i} (x) -\nu_{\gamma}(f))- \sum_{i=1}^k Z_i (x,y) \right| = o (n^{\max(\gamma,1/4)} (\log n)^{\eta + \varepsilon}) \text{ $\nu_{\gamma}\otimes \lambda$ a.e.}
\end{equation}
where $\eta =3/4$ for $\gamma \leq 1/4$, $\eta =5/6$ if  $\gamma =1/3$ and $\eta =1/2$ otherwise.
\item For any $\gamma \in ] 0, 1/2[$, one can construct on  the probability space $([0,1]\times[0,1], \nu_{\gamma}\otimes \lambda)$ a sequence of i.i.d.  Gaussian random variables $(Z_i)_{i \geq 1}$ centered and with variance $\sigma^2$, such that
\begin{equation}\label{L2}
\left(\int \sup_{k \leq n} \left |\sum_{i=1}^k (f\circ T_{\gamma}^{i} (x) -\nu_{\gamma}(f))- \sum_{i=1}^k Z_i (x,y) \right|^2 \! \! \! 
\nu_{\gamma}(dx) dy  \! \right) ^{\frac 1 2} \! \! \! = o (n^{\max(\gamma,1/4)} (\log n)^{\eta }) \, ,
\end{equation}
where  $\eta =1/4$ if $\gamma <1/4$, $\eta =1/2$ if  $\gamma = 1/4$ or $\gamma =1/3$ and $\eta =1/2-\gamma$ otherwise.
\end{itemize}

\subsection{Maps modelled by Young towers}\label{tau}  Let $(M,d)$ be a complete bounded separable metric space with the Borel $\sigma$-algebra.
Let $T: M \mapsto M$ be a map that can be modelled by a Young tower (see \cite{Young99}, or for instance \cite{CDKM20}), and denote by $\nu$ the  $T$-invariant probability measure on $M$ induced by this Young tower. 
Let $\varphi$ be an H\"older observable from $M$ to ${\mathbb R}$. Then, as  explained  in \cite{DMR22}, one can compute the coefficients $\theta(k)$ of the sequence $(\varphi (X_i))_{i \geq 0}$, where $(X_i)_{i \geq 0}$ is the stationary Markov chain associated with $T$, whose transition kernel is the Perron-Frobenius operator of the composition by $T$ with respect to $\nu$. We infer that Theorems \ref{thSAPfaible}  and \ref{thSAP} apply to   $(\varphi (X_i)-\nu(\varphi))_{i \geq 0}$, and also (proceeding as in \cite[Section 3]{MR12}) to the sequence    $(\varphi (T^i)-\nu(\varphi))_{i \geq 0}$ on the probability space $(M, \nu)$. 

More precisely, the behaviour of the coefficients $\theta(k)$ of the sequence $(\varphi (X_i))_{i \geq 0}$ depends on the behaviour of the  return time to the base of the tower. For instance, if the return time has a moment of order $p>1$, then $\sum_{k>0} k^{p-2} \theta(k) < \infty$;   if it has a weak moment of order $p>1$  then  $\theta(k) = O(k^{1-p})$; if it has an exponential moment, then  $\theta(k)=O(a^k) $, for some $a \in ]0,1[$ (see  Proposition 2.1 and its remark 2.1 in \cite{CDM23} in case of strong moments and Proposition 5.3 in \cite{DM} in case of weak moments). Our results apply in particular to H\"older observables of the LSV map, leading to the same upper bounds \eqref{AS}-\eqref{L2} as for BV observables. 

Note that, for H\"older observables of maps that can be modelled by a Young tower, optimal almost sure rates can be obtained via another method described in the paper \cite{CDKM20} and inspired by \cite{BLW}. In particular, for H\"older observables of the LSV map of parameter $\gamma < 1/2$, the optimal  almost sure rate $o(n^\gamma (\log n)^{\gamma+ \varepsilon})$ is given in \cite{CDKM20}.  However, no results similar to the ${\mathbb L}^2$ control \eqref{L2} are given in \cite{CDKM20}. Moreover, with our approach, we can also obtain rates for a larger class of continuous observables (including H\"older observables of any exponent) when the return time to the base as an exponential moment (see \cite[Corollary 3.4]{DMR22}). 

\subsection{Rates of convergence in the functional CLT}  \label{SectionRFCLT}
Let $(X_i)_{i \in {\mathbb Z}}$ be a strictly stationary sequence of centered and bounded random variables, and define the Donsker line
\[
B_n  (t) = \frac{1}{\sqrt{n}} \Big (  \sum_{k=1}^{[nt]} X_k  + (nt - [nt]) X_{[nt]} \Big )  \, .
\]
It is well known that, if $\sum_{k>0} \theta(k) < \infty$, then $B_n $ converges in distribution to $\sigma B$ on the space $C([0,1])$ with the uniform distance, where $B$ is a standard Brownian motion, and $\sigma^2$ is the covariance series defined in Theorem  \ref{Fuk-inequality-theta}.

One can now ask for the convergence rate in this functional CLT, with respect to the Wasserstein distance of order $p\geq 1$, that is rates for the quantity
$
   W_p(P_{B_n }, P_{\sigma B})
$, 
where $P_{B_n }$ and $P_{\sigma B}$ are the distributions of $B_n $ and $\sigma B$, and the cost function is 
$|\cdot|^p_\infty$, $|\cdot |_\infty$ being the supremum norm on $C([0,1])$. Note that, by definition of $W_p$,
\[
W_{p} ( P_{B_n } , P_{\sigma B} )   \leq  \Big \Vert  \sup_{0 \leq t \leq 1}  | B_n ( t) - \sigma B(t) | \Big \Vert_{p} \, ,
\]
for any standard Brownian motion $B$.
Consequently Theorem  \ref{thSAPfaible} applies when $p=2$. For instance, if  $\sum_{k >0} k^2 \theta(k) < \infty$,  Item (i) of Theorem \ref{thSAPfaible} implies that  
\begin{equation} \label{functional}
W_2  ( P_{B_n } , P_{\sigma B} )   = O   \bigl( n^{-1/4}  (\log n)^{1/4} \bigr)  \, .
\end{equation}

Note that, according to Section \ref{tau}, the upper bound \eqref{functional} applies to sequences $(\varphi (T^i)-\nu(\varphi))_{i \geq 0}$, where $\varphi$ is an H\"older observable and $T$ can be modelled by a Young tower with return time to the base having a moment of order 4. 
Let us compare this result with some recent results obtained in this context by Liu and Wang  \cite{LW22}. 

  Assume that the return time to the base  has a finite  moment of order $q \geq 4$. Using the Skorohod embedding theorem, Liu and Wang  \cite{LW22} (see their Theorem 3.4)  proved  that 
\beq \label{Wp2LW}
W_{q/2} ( P_{B_n } , P_{\sigma B} )    = O ( n^{-(q-2)/ (4(q-1))}) \, .
\eeq
Therefore, for $q=4$,  \eqref{Wp2LW}  gives the rate $O \bigl(n^{-1/6} \bigr)$, while the upper bound  \eqref{functional} gives the rate  
$O   \bigl( n^{-1/4}  (\log n)^{1/4} \bigr)$.

\section{Proof of Theorem \ref{Fuk-inequality-theta}} \label{SectProofFN}

\setcounter{equation}{0}

Starting from inequality (A.42) in \cite{MR12} (in the bounded case) together with the fact that $\Vert \BBE_{i-q} (X_i) \Vert_1 \leq \theta(q) $, we infer that for any nondecreasing, non negative and convex function $\varphi $ and 
any $x >0$, 
\[ 
\BBP \big ( S_n^* \geq 4x \big ) \leq  \frac{\BBE (\varphi ( S_n ))}{\varphi (x)} + n x^{-1} \theta ([x]) \, .
\]
Now, since $x^3  \leq  8 \sum_{k=1}^{[x]} k^2$ when $ x\geq 1$,  $(\theta(k))_{k\geq 0}$ is non increasing and $\theta(0) \leq 1$ (since 
$\Vert X_0 \Vert_{\infty} =1$), 
\[
x^3 \theta ([x]) \leq 1 + 8 \sum_{k=1}^{[x]} k^2 \theta (k) \leq 1 + 8 \sum_{k \geq 1} k ( k \wedge x)  \theta (k) \, . 
\] 
Hence
\beq \label{1F}
\BBP \big ( S_n^* \geq 4x \big ) \leq  \frac{\BBE (\varphi ( S_n ))}{\varphi (x)} + \frac{8n}{x^4} 
\Bigl( 1 + \sum_{k \geq 1} k ( k \wedge x)  \theta (k) \, \Bigr) \, . 
\eeq
\par
Next we handle the first term in the right-hand side of \eqref{1F} with the following selection of $\varphi$: For any 
real $t$, 
\[
\varphi(t) =  \left\{
\begin{array}{ll}0 & \text{if $ t \leq x/2$} \\
\frac{1}{24} ( t - \frac{x}{2} )^4 & \text{if $x/2 \leq t \leq x $} \\
\frac{x^4}{24 \times 2^4} + \frac{x}{12} ( t - x)^3+  \frac{x^2}{16} ( t - x)^2 + \frac{x^3}{48} ( t -x) & \text{if $t \geq x $} \, .
\end{array}
\right.
\]
This is a nondecreasing and convex function such that 
$ \|\varphi^{(3)}\|_{\infty}=x/2$ and  $\|\varphi^{(4)}\|_{\infty} = 1$. Furthermore $384 \, \varphi (x) =  x^4$, whence 
$\BBE (\varphi ( S_n ))/\varphi (x) = 384\, x^{-4} \BBE (\varphi ( S_n ))$. 
\par\ssk
To bound up $ \BBE (\varphi ( S_n ))$, we start by a symmetrization argument. Let $(X'_k)_{k\in {\mathbb Z}}$ be a stationary sequence 
independent of the sequence $\XZ$ and with the same joint law as $\XZ$.
Set $S'_n = X'_1 + X'_2 + \cdots + X'_n$. Since $S'_n$ is centered and independent of $S_n$, it follows from
the conditional version of the Jensen inequality that 
\beq \label{Symm}
\BBE \bigl(\varphi (S_n)\bigr) \leq \BBE  \bigl(\varphi (S_n- S'_n)\bigr) \, .
\eeq
Hence Inequality \eqref{ineFN1} will follow from \eqref{1F} if we  prove that 
\beq \label{2F}
\BBE \bigl(\varphi (S_n - S'_n)\bigr) \ll   {\bf 1}_{ \{\sigma^2 >0 \} }  \Big (  \frac{n \sigma^2}{x} \Big )^4  
\exp \Big (- \frac{x^2}{ 16 n \sigma^2 } \Big )    + n \Bigl( \Theta_1 \Theta_2  +  
\sum_{k \geq 1}  k (k \wedge x )  \theta(k)  \Bigr) \, . 
  \eeq
\par
Define then the stationary sequence $\ZZ$ of centered and bounded random variables by 
$Z_i = X_i - X'_i$ for any integer $i$, and set 
\beq \label{deftildeS}
\tilde S_0 = 0 \ \text{ and } \tilde S_n = Z_1+Z_2 +  \cdots + Z_n = S_n - S'_n \text{ for any integer } n>0.
\eeq 
From the definition of $\ZZ$,  
\beq \label{CovZ}
\BBE ( Z_0 Z_k) = 2 \BBE ( X_0X_k) \text{ for any } k\in \BBN, \text{ whence } \BBE ( Z_0^2) + 
2 \sum_{k\geq 1} \BBE ( Z_0 Z_k ) = 2 \sigma^2 .
\eeq
\par
To prove \eqref{2F}, we shall apply the Lindeberg method to $\ZZ$: we 
consider a sequence $(Y_k)_{k \geq 1}$ of i.i.d.  random variables with normal law $N (0, 2\sigma^2)$,  independent of 
$\ZZ$. Set  $T_0=0$ and $T_n = Y_1 + Y_2 + \cdots + Y_n$ for $n>0$. Clearly 
\beq \label{FromStoT}
\BBE ( \varphi ( \tilde S_n ) ) = \BBE ( \varphi (T_n) )  + \bigl( \BBE ( \varphi ( \tilde S_n ) ) - \BBE ( \varphi (T_n) )  \bigr).
\eeq
We start by computing  $ \BBE ( \varphi (T_n) ) $. If $\sigma^2 = 0$, then $T_n=0$ and 
$ \BBE ( \varphi (T_n) ) = 0$. If $\sigma^2>0$, then $T_n$ has the normal law $N (0, 2n\sigma^2 )$ and consequently 
\[
\BBP \bigl(   T_n  \geq t  +x/2  \bigr) \leq  \exp \bigl( -  (t+x/2)^2 /  (4n\sigma^2)   \, \bigr) \, ,
\]
Since $(t+x/2)^2 \geq tx  + x^2/4$, we derive that  
\[
 \int_0^{+ \infty} t^3   \BBP \bigl(   T_n  \geq t  +x/2  \bigr)   dt  \leq e^{-x^2/(16n\sigma^2 )}  \int_0^{+ \infty} t^3   e^{-  tx / (4n\sigma^2 )} dt   \leq 6  \Big (  \frac{4n \sigma^2}{ x}\Big )^4  e^{-x^2/(16n\sigma^2 )}     \, .
\]
So, overall,  
\[
\BBE ( \varphi (T_n) )  \leq  \BBE ( (T_n -x/2)^4_+ )  \leq   4 \int_0^{+ \infty} t^3 \BBP \big ( ( T_n  \geq t +x/2  \big )  dt    \ll  \Big ( \frac{n  \sigma^2 }{ x} \Big )^{4}  e^{-x^2/( 16 n \sigma^2 )} \, .
\]
\par
According to \eqref{FromStoT} and the above inequality, to end the proof of  Inequality \eqref{2F}, 
it remains to prove that 
\begin{equation} \label{aimFuk}
\BBE ( \varphi (\tilde S_n) ) - \BBE ( \varphi (T_n) )  \ll n 
\bigl (  \Theta_1 \Theta_2 + \sum_{k \geq 1} k ( k \wedge x) \theta(k)  \bigr )    \, .
\end{equation} 
With this aim, let $\varphi_{n-k} (t) = {\mathbb E} ( \varphi (t  + T_n - T_k)  )$ and define, for $k \geq 1$, 
\[
\Delta_{n,k} =  \varphi_{n-k} ( \tilde S_{k-1} + Z_k) - \varphi_{n-k} ( \tilde S_{k-1} + Y_k) \, .
\]
The functions $\varphi_{n-k} ( \cdot)$  are ${\mathcal C}^{\infty}$, $\Vert  \varphi^{(3)}_{n-k} \Vert_{\infty}  =b_3 \leq x/2$ and $\Vert  \varphi^{(4)}_{n-k} \Vert_{\infty} =b_4 \leq 1$. 
Since the sequence $\Y$ is independent of the sequence $\ZZ$,
\begin{equation} \label{sumdelta-F}
\BBE ( \varphi (\tilde S_n) ) - \BBE ( \varphi (T_n) )= \sum_{k=1}^n {\mathbb E} ( \Delta_{n,k} ) .
\end{equation}
\par
\begin{Notation} Set  $\Delta_{n,k}^{(1)} =  \varphi'_{n-k} (\tilde S_{k-1}) (Z_k-Y_k)$, 
$\Delta_{n,k}^{(2)} = \varphi''_{n-k} (\tilde S_{k-1} ) (Z_k^2 - Y_k^2)$  
and $\Delta_{n,k}^{(3)} = \varphi_{n-k}^{(3)} (\tilde S_{k-1} ) (Z_k^3 - Y_k^3) $.  Let 
$\Delta^*_{n,k} = \Delta_{n,k}^{(1)} +  \frac{1}{2}    \Delta_{n,k}^{(2)} +   \frac{1}{6}  
\Delta_{n,k}^{(3)}$.
\end{Notation}
With the above notations, from the Taylor integral formula at order $4$, 
\beq \label{deltaetoileF}
 \Delta_{n,k} = \Delta^*_{n,k} + R_{n,k} \, ,
\eeq
\[\text{with }
R_{n,k} =\frac{1}{6}    Z_k^4 \int_0^1 (1-s)^3 \varphi_{n-k}^{(4)} (\tilde S_{k-1}  + s Z_k) ds  
+   \frac{1}{6}   Y_k^4  \int_0^1 (1-s)^3 \varphi_{n-k}^{(4)} (\tilde S_{k-1}  + s Y_k) ds  \, .
\]
Since $0 \leq \varphi_{n-k}^{(4)}  \leq 1$, $\Vert Z_k \Vert_{\infty}  \leq 2 $ and 
$\BBE ( Y_k^4) = 12 \sigma^4$, we derive that 
\begin{equation}\label{reste1-F}
0 \leq \sum_{k=1}^n \BBE ( R_{n,k} )  \leq n ( 1 + \sigma^4)  \leq 4 n ( 1 + \Theta_1 \Theta_2)   \, .
\end{equation}
Indeed  $\sigma^4 \leq  4 \Theta_1^2 \leq 4 \Theta_1 \Theta_2$.   Next
\begin{equation} \label{linddec1-F}
\BBE (   \Delta^*_{n,k}   ) = \BBE (   \Delta_{n,k}  -  R_{n,k}  ) =  \BBE ( \Delta_{n,k}^{(1)}) + {\textstyle \frac{1}{2} }   \BBE ( \Delta_{n,k}^{(2)}) +  {\textstyle \frac{1}{6} }   \BBE ( \Delta_{n,k}^{(3)}) \, .
\end{equation}
From the fact that $(Y_k)_{k\geq 1}$ is independent of $\ZZ$, 
$\BBE ( \Delta_{n,k}^{(1)}) = \BBE ( \varphi'_{n-k} (\tilde S_{k-1}) Z_k)$, 
$\BBE ( \Delta_{n,k}^{(2)}) = \BBE ( \varphi''_{n-k} (\tilde S_{k-1}) (Z_k^2 - 2 \sigma^2) )$ and 
 $\BBE ( \Delta_{n,k}^{(3)}) = \BBE ( \varphi_{n-k}^{(3)} (\tilde S_{k-1} ) Z_k^3 )$. 
We now develop each term in the right-hand side of \eqref{linddec1-F} with the help of the Lindeberg method. 
\par\ssk
Since $\BBE ( \varphi'_{n-k} (0) Z_k  )= \varphi'_{n-k} (0)  \BBE (Z_k) =0$,  we have 
\begin{align} \label{linddec2-F}
 \BBE ( \Delta_{n,k}^{(1)}) & =  \sum_{i=1}^{k-1}  \BBE \big ( \{  \varphi'_{n-k} (\tilde S_{k-i}) - \varphi'_{n-k} (\tilde S_{k-i- 1})  \} Z_k   \big )   \nonumber \\
&  =    \BBE  (\Delta_{n,k,2}^{(1)} )  +  {\textstyle \frac{1}{2}}   \BBE  (\Delta_{n,k,3}^{(1)} )     +  A_{n,k,2}^{(1)}  +   
{\textstyle \frac{1}{2}}   A_{n,k,3}^{(1)} + B_{n,k}^{(1)}     \, ,
\end{align}
where the following notations have been used:  for $j= 2,3$, 
\begin{align*} 
 \Delta_{n,k,j}^{(1)} & =  \sum_{i=1}^{k-1}   \varphi^{(j)}_{n-k} (\tilde S_{k-i- 1})     ( Z_{k-i}^{j-1}Z_k )^{(0)},    \
A_{n,k,j}^{(1)} = \sum_{i=1}^{k-1}    \BBE \{  \varphi^{(j)}_{n-k} (\tilde S_{k-i- 1})  \}    \BBE( Z_{k-i}^{j-1}Z_k ) ,
\cr &
\text{ and } 
B_{n,k}^{(1)}   = 
\frac{1}{2}  \sum_{i=1}^{k-1}  \int_0^1 (1-s)^2 \BBE \big (   \varphi_{n-k}^{(4)} (\tilde S_{k-i-1}  + s Z_{k-i})  Z_{k-i}^3 Z_k  \big ) ds \, . \cr
\end{align*}
In the  decomposition \eqref{linddec2-F}, $A_{n,k,3}^{(1)} = 0$. This is due to the lemma below, 
whose proof uses the symmetry properties of $\ZZ$. 

\begin{Lemma} \label{ThirdMomentsZ} For any integers $i$, $j$ and $k$, 
$\BBE (Z_i Z_j Z_k) = 0$.
\end{Lemma}
\noindent
{\bf Proof of Lemma \ref{ThirdMomentsZ}.} From the definition of $\ZZ$, $(-Z_i , -Z_j, -Z_k)$
has the same joint law as $(Z_i, Z_j, Z_k)$. Hence $\BBE (-Z_i Z_j Z_k) = \BBE ( Z_i Z_jZ_k)$, which implies 
Lemma \ref{ThirdMomentsZ}. \qed
\par\medskip
Next, concerning the coefficients $(\theta_{Z,p,q} (k))_{k \geq 0}$ associated with the sequence $\ZZ$, they can be compared with the former coefficients $(\theta_{X,p,q} (k))_{k \geq 0}$ as follows.
\begin{Lemma} \label{coeffZ} For any integer $k\geq 0$ and any positive integers $p$ and $q$, 
\[
\theta_{Z,p,q} (k) \leq 2^{q+1} \theta_{X,p,q} (k) \, .
\]
\end{Lemma}
\noindent
{\bf Proof of Lemma \ref{coeffZ}.}  Let ${\mathcal F}'_0 = \sigma ( X_i', i \leq 0 )$ and  $ {{\tilde {\mathcal  F}_0 }} = {\mathcal F}_0 \vee {\mathcal F}'_0$. By definition, 
\[
\theta_{Z,p,q} (k)  = \sup_{ k_p>k_{p-1}> \ldots > k_2>k_1 \geq k \atop (a_1, \dots, a_p ) \in \Gamma_{p,q}} \Big \Vert  \BBE \Big ( \prod_{i=1}^p Z_{k_i}^{a_i} |   {{\tilde {\mathcal  F}_0 }}  \Big  ) - \BBE  \Big ( \prod_{i=1}^p Z_{k_i}^{a_i}  \Big )  \Big \Vert_1 \, .
\]
Let  $(a_i)_{ 1 \leq i \leq p} $ in $ \Gamma_{p,q}$ and $k_1, \cdots, k_p$ such that $k_p>k_{p-1}> \ldots > k_2>k_1 \geq k$. Note  that 
$Z_{k_i}^{a_i} = \sum_{\ell =0}^{a_i} C_{a_i}^{\ell} X_{k_i}^{\ell}  (-X'_{k_i}) ^{a_i-\ell}  $. Next,  let  $(\alpha_i)_{ 1 \leq i \leq p} $, $(\beta_i)_{ 1 \leq i \leq p} $ in $ \Gamma_{p,q}$ and $k_1, \cdots, k_p$ such that $k_p>k_{p-1}> \ldots > k_2>k_1 \geq k$. Setting $U =\prod_{i=1}^p X_{k_i}^{\alpha_i}$,  $V =\prod_{i=1}^p (X'_{k_i})^{\beta_i}$ and
$\sigma (X') = \sigma (X'_k, k \in {\mathbb Z} )$, note that 
\begin{align*}
 \big \Vert  \BBE \big (UV  |   {{\tilde {\mathcal  F}_0 }}  \big  ) - \BBE  \big ( UV \big )  \big \Vert_1 
 & =   \big \Vert  \BBE \big (  \BBE \big ( UV  |  \sigma (X') \vee  {{{\mathcal  F}_0 }}  \big  )   |   {{\tilde {\mathcal  F}_0 }}  \big  ) - \BBE  \big ( UV \big )  \big \Vert_1  \\
&  =    \big \Vert  \BBE \big (  V  \BBE \big ( U  |  \sigma (X') \vee  {{{\mathcal  F}_0 }}  \big  )   |   {{\tilde {\mathcal  F}_0 }}  \big  ) - \BBE  (V) \BBE  (U)   \big \Vert_1 \, .
\end{align*}
Next we use the following well-known fact. Let $Y$ be  an integrable random variable, and ${\mathcal G}_1$ and ${\mathcal G}_2$ be two $\sigma$-algebras such that $\sigma(Y) \vee {\mathcal G}_1$ is independent of ${\mathcal G}_2$, then 
\[
\BBE \big ( Y  |  {\mathcal G}_1 \vee {\mathcal G}_2  \big  )  = \BBE \big ( Y  |  {\mathcal G}_1 \big  )  \quad \text{a.s.}
\]
From the above fact and since we assume that the $X_i$'s are uniformly bounded by one, 
\begin{align*}
 \big \Vert  \BBE \big (UV  |   {{\tilde {\mathcal  F}_0 }}  \big  ) - \BBE  \big ( UV \big )  \big \Vert_1 
&  =    \big \Vert  \BBE \big (  V   \big \{  \BBE \big ( U  |    {{{\mathcal  F}_0 }}  \big  )   - \BBE (U) \big \}  |   {{\tilde {\mathcal  F}_0 }}  \big  ) 
+    \big \{   \BBE \big (  V   |   {{\tilde {\mathcal  F}_0 }}  \big  )   - \BBE  (V)  \big \} \BBE  (U)   \big \Vert_1  \\
& \leq   \big \Vert  \BBE \big ( U  |    {{{\mathcal  F}_0 }}  \big  )   - \BBE (U)  \big \Vert_1
+    \big \Vert   \BBE \big (  V   |   {\mathcal  F}'_0   \big  )   - \BBE  (V)    \big \Vert_1 \leq 2  \theta_{X,p,q} (k)  \, .
\end{align*}
So, overall,
\[
\theta_{Z,p,q} (k)  \leq 2  \sup_{ (a_1, \dots, a_p ) \in \Gamma_{p,q}}  \prod_{i=1}^p 2^{a_i}   \theta_{X,p,q} (k)  \, , 
\]
proving the lemma.  \qed
\par\medskip
 We now handle the rests $B_{n,k}^{(1)}$.  
 Since $\Vert Z_{k-i}^3  \varphi_{n-k}^{(4)} (\tilde S_{k-i-1}  + s Z_{k-i})   \Vert_{\infty} \leq 2 $, 
\beq \label{BoundB1-F}
\sum_{k=1}^n  \vert B_{n,k}^{(1)}  \vert \ll  n   \Theta_1   \, .
\eeq
To handle the second term in the right-hand side of  \eqref{linddec1-F}, we introduce  the following additional notations.

\begin{Notation} 
Let $ \gamma_i= \BBE (Z_0Z_i) $.  Define 
$\beta_{2, {k}} = 2 \sum_{i =1}^{k-1} \gamma_i$ and $ \beta_{2}^{({k})} = 2 \sum_{i \geq {k} }  \gamma_i$.  
\end{Notation}

Since  $\BBE( Z_{k-i}Z_k ) = \gamma_i$, note that 
\begin{align} \label{Boundr21-F-beta2}
\frac{1}{2}     \BBE \{  \varphi''_{n-k} (\tilde S_{k-1})  \}   \beta_{2,k}   - A_{n,k,2}^{(1)}   
& =    \sum_{i=1}^{k-1}  \gamma_i   \sum_{j=1}^i \BBE \{  \varphi''_{n-k} (\tilde S_{k-j})  -  \varphi''_{n-k}  (\tilde S_{k-j -1})   \}      \nonumber \\
& =    \sum_{i=1}^{k-1}   \gamma_i  \sum_{j=1}^i \BBE  \big \{   \varphi^{(3)}_{n-k}  (\tilde S_{k-j -1})   Z_{k-j}  \big \}  + r_{n,k,2}^{(1)}  \, , \end{align}
where 
\[
 r_{n,k,2}^{(1)} :=   \int_0^1 (1-t)  \sum_{i=1}^{k-1} \gamma_i   \sum_{j=1}^i \BBE  \big \{    \varphi^{(4)}_{n-k}  (\tilde S_{k-j -1} + t Z_{k-j}) Z_{k-j}^2  \big \}     dt \, .
\]
From the fact that $ \Vert   \varphi^{(4)}_{n-k}  (\tilde S_{k-j -1} + t Z_{k-j}) Z_{k-j}^2 \Vert_1 \leq b_4 \BBE (Z_0^2) \leq 4 b_4$, it follows that 
\beq \label{Boundr21-F}
\sum_{k=1}^n \vert  r_{n,k,2}^{(1)} \vert  \ll  n  \sum_{i=1}^{n} i \gamma_i \ll n  \Theta_2 . 
\eeq
Starting from \eqref{linddec1-F} and taking into account \eqref{linddec2-F}, \eqref{BoundB1-F}, \eqref{Boundr21-F} and the definition of $\sigma^2$, we get 
\begin{multline} \label{linddec3-F}
\BBE (   \Delta^*_{n,k} )   
 =  \BBE  (  \Delta_{n,k,2}^{(1)}  ) + { \frac{1}{2} }  \BBE \big (   \varphi''_{n-k} (\tilde S_{k-1} ) (Z_k^2 )^{(0)}\big )  -  { \frac{1}{2} }  \BBE \big (   \varphi''_{n-k} (\tilde S_{k-1} ) \big )  \beta_2^{({k})} 
 + \frac{1}{2}   \BBE  (\Delta_{n,k,3}^{(1)} )  
 \\
   -  \sum_{i=1}^{k-1}  \gamma_i  \sum_{j=1}^i \BBE  \big \{   \varphi^{(3)}_{n-k}  (\tilde S_{k-j -1})   Z_{k-j}  \big \}   + 
\frac{1}{6}   \BBE \big (  \varphi_{n-k}^{(3)} (\tilde S_{k-1} )  Z_k^3  \big )  
 + \Gamma^{(1)}_{n,k} \, ,
\end{multline}
where $  \Gamma^{(1)}_{n,k}   $ satisfies
\beq \label{boundGamma-F}
\sum_{k=1}^n \vert  \Gamma^{(1)}_{n,k}   \vert  \ll n   \Theta_2 \, .
\eeq
Note that 
\[
\big |  \BBE \big (   \varphi''_{n-k} (\tilde S_{k-1} ) \big )\big | 
  \leq  \BBE \big (  \tilde S_{k-1}  + T_n - T_k \big )^2 =  2 \BBE (S^2_{k-1} ) + 2(n-k) \sigma^2  
  \leq 4 n   \Theta_1  \, .
\]
Hence
\beq \label{B0term1-F}
\sum_{k=1}^n   \vert \BBE \big (   \varphi''_{n-k} (\tilde S_{k-1} ) \big )  \beta_2^{({k})}  \vert  \ll     n    \Theta_1  \sum_{k=1}^n  \sum_{i \geq {k}  }   \theta  (i)   \ll  
n     \Theta_1 \Theta_2 \,  .
\eeq
We now handle the quantity $  \BBE \big ( \varphi''_{n-k} (\tilde S_{k-i-1})  (Z_{k-i}Z_k  )^{(0)} \big ) $ appearing in the first  two terms of the right hand side of \eqref{linddec3-F} for any integer  $i$  in $[0, k-1]$. With this aim, noticing that $ \BBE \big ( \varphi''_{n-k} (0)  (Z_{k-i}Z_k  )^{(0)} \big ) =0$, 
$$
\BBE   \big \{   \varphi''_{n-k} (\tilde S_{k-i- 1})    ( Z_{k-i}Z_k )^{(0)}   \big \} =  \sum_{j=i+1}^{k-1 } 
\BBE   \big \{   (  \varphi''_{n-k} (\tilde S_{k-j }) -  \varphi''_{n-k} (\tilde S_{k-j -1})  )   ( Z_{k-i}Z_k )^{(0)}   \big \} ,
$$
whence
\begin{multline}  \label{linddec4-F} 
\BBE   \big \{   \varphi''_{n-k} (\tilde S_{k-i- 1})    ( Z_{k-i}Z_k )^{(0)}   \big \}
= \sum_{j=i+1}^{k-1} \BBE   \big \{     \varphi^{(3)}_{n-k}(\tilde S_{k-j -1})   Z_{k-j}  ( Z_{k-i}Z_k )^{(0)}   \big \} \\
+  \sum_{j=i+1}^{k-1 } \int_0^1 (1-t)  \BBE   \big \{      \varphi^{(4)}_{n-k}(\tilde S_{k-j -1} + t Z_{k-j})   Z^2_{k-j}  ( Z_{k-i}Z_k )^{(0)}   \big \} dt  \, .
\end{multline}
Since $ \Vert  \varphi^{(4)}_{n-k}(\tilde S_{k-j -1} + t Z_{k-j})     \Vert_{\infty} \leq 1$, it follows that
\begin{multline} \label{B2term1-F}
\sum_{k=1}^n \sum_{i=0}^{{k-1}}   \sum_{j=i+1}^{k-1 } \vert \BBE   \big \{     \varphi^{(4)}_{n-k}(\tilde S_{k-j -1} + t Z_{k-j})   Z^2_{k-j}  ( Z_{k-i}Z_k )^{(0)}    \big \}  \vert \\
  \ll \sum_{k=1}^n \sum_{i=0}^{{k-1}}   \sum_{j=i+1}^{k-1 } 
\theta (j-i) \wedge \theta (i)  
\ll  n  \Theta_2   \, .
\end{multline}
Next, set 
$ \Delta_{n,k,2}^{(1,3)}  (i) = \sum_{j=i+1}^{k-1 }     \varphi^{(3)}_{n-k}(\tilde S_{k-j -1})   Z_{k-j}  ( Z_{k-i}Z_k )^{(0)}$.
Starting from \eqref{linddec3-F} and taking into account \eqref{boundGamma-F},  \eqref{B0term1-F}, \eqref{linddec4-F} and \eqref{B2term1-F}, we get
\begin{multline} \label{linddec5-F}
\BBE (   \Delta^*_{n,k} )   
 =  { \frac{1}{2} }   \BBE (   \Delta_{n,k,2}^{(1,3)}  (0) )  + \sum_{i=1}^{{k-1}} \BBE (    \Delta_{n,k,2}^{(1,3)}  (i)  )   +
    \frac{1}{2}   \BBE  (\Delta_{n,k,3}^{(1)} )   \\
    -  \sum_{i=1}^{k-1}  \gamma_i  \sum_{j=1}^i \BBE  \bigl( \varphi^{(3)}_{n-k}  (\tilde S_{k-j -1})   Z_{k-j}  \bigr)   
     + 
\frac{1}{6}   \BBE \big (  \varphi_{n-k}^{(3)} (\tilde S_{k-1} )  Z_k^3   \big )  
+ \Gamma^{(2)}_{n,k} \, ,  
\end{multline}
where the rests $\Gamma^{(2)}_{n,k}$ satisfy
\beq  \label{linddec5-F-rests}
\sum_{k=1}^n \vert  \Gamma^{(2)}_{n,k}   \vert  \ll n   \Theta_2 +  n  \Theta_1  \Theta_2     \, .
\eeq
Introduce now the following notations:  for  any integer $i $, let 
\[
{\tilde  \Delta}_{n,k,2}^{(1,3)}  (i) := \sum_{j=i+1}^{k-1 }   \big \{     \varphi^{(3)}_{n-k}(\tilde S_{k-j -1})   (  Z_{k-j}  ( Z_{k-i}Z_k )^{(0)}  )^{(0)}  \big \}  \, .
\]
Starting from \eqref{linddec5-F},  and using the fact that $\BBE    \big \{ Z_{k-j}   Z_{k-i}Z_k   \big \}  =  
\BBE    \big \{ Z_{k-j}  ( Z_{k-i}Z_k )^{(0)}   \big \}  = 0$ thanks to Lemma \ref{ThirdMomentsZ}, we get that
\begin{multline} \label{linddec6-F}
\BBE (   \Delta^*_{n,k} )   
 =  { \frac{1}{2} }  \BBE ({\tilde  \Delta}_{n,k,2}^{(1,3)}  (0) ) + \sum_{i=1}^{{k-1}}   \BBE ({\tilde  \Delta}_{n,k,2}^{(1,3)}  (i))   +  
\frac{1}{6}   \BBE \big (  \varphi_{n-k}^{(3)} (\tilde S_{k-1} )  Z_k^3 ) \big )   \\
  +    \frac{1}{2}    \BBE  (\Delta_{n,k,3}^{(1)} )  
    -  \sum_{i=1}^{k-1}  \gamma_i  \sum_{j=1}^i \BBE  \big \{   \varphi^{(3)}_{n-k}  (\tilde S_{k-j -1})   Z_{k-j}  \big \}    + \Gamma^{(2)}_{n,k} . 
\end{multline}
We handle now the first two terms in the right-hand side of \eqref{linddec6-F}. Since $\Vert   \varphi^{(3)}_{n-k} \Vert_{\infty} \leq x/2$,
\[
 \big | \BBE \big \{  \varphi^{(3)}_{n-k}(\tilde S_{k-j -1 })   Z_{k-j}  ( Z_{k-i}Z_k )^{(0)}  \big \} \big |  \ll x  \{ \theta(j-i) \wedge \theta(i) \} \, ,
\]
whence
\beq \label{B5dec6F1}
 |    \BBE ( {\tilde  \Delta}_{n,k,2}^{(1,3)}  (i)  )  |  \ll  x  \Big ( i  \theta(i)  + \sum_{j =i }^{k-i}  \theta(j) \Big )    \, .
\eeq
On another hand, since  $ \BBE   \big \{   \varphi^{(3)}_{n-k} (0)    Z_{k-j}  ( Z_{k-i}Z_k )^{(0)}    \big \}   =0$, we write
  \begin{multline*}
\BBE \big \{  \varphi^{(3)}_{n-k}(\tilde S_{k-j -1 })    Z_{k-j}  ( Z_{k-i}Z_k )^{(0)}  \big \} \\
= \sum_{u=1}^{k-j-1}  \BBE \big \{  \big ( \varphi^{(3)}_{n-k}(\tilde S_{k-j -u })  -     \varphi^{(3)}_{n-k}( \tilde S_{k-j -u -1 })   \big )  Z_{k-j}  ( Z_{k-i}Z_k )^{(0)}   \big \} \\
=  \sum_{u=1}^{k-j-1}   \int_0^1 \BBE \big \{   \varphi^{(4)}_{n-k}(\tilde S_{k-j -u -1 }   +t  Z_{k-j -u  }  )   Z_{k-j -u }  
Z_{k-j}  ( Z_{k-i}Z_k )^{(0)}   \big \} dt   \, .
\end{multline*}
Since $\Vert   \varphi^{(4)}_{n-k}(S_{k-j -u -1 }   +t  X_{k-j -u  }  )  \Vert_{\infty}  \leq 1$,  it follows that 
\begin{equation} \label{B5dec6F2}
|    \BBE ( {\tilde  \Delta}_{n,k,2}^{(1,3)}  (i)  )  |
% \\ \leq   \sum_{i = 0}^{k}  \sum_{j=i+1}^{k-1}  \sum_{u=1}^{k-i-1}    \int_0^1 \big  | \BBE \big \{   \varphi^{(4)}_{n-k}(S_{k-j -u -1 }   +t  X_{k-j -u  }  )   X_{k-j -u }  (  X_{k-j}  ( X_{k-i}X_k )^{(0)}  )^{(0)}  \big \}  \big | dt  \\
 \ll      \sum_{j=1}^{k }  \sum_{u=1}^{k}  \big ( \theta(u) \wedge  \theta(j)    \wedge  \theta(i)  \big )  \ll     \sum_{j = 1}^k   j  \big ( \theta(j) \wedge  \theta(i)  \big )     \, .
\end{equation}
Therefore, using  the upper bound  \eqref{B5dec6F1} when $i \geq x$ and the upper bound \eqref{B5dec6F2} when $ i < x$, we derive
 \begin{multline*}
 \sum_{i=0}^{{n}} |    \BBE ( {\tilde  \Delta}_{2,3,i}^{(1)}  )  | \ll x  \sum_{i \geq [x]} i  \theta(i)   + 
 \sum_{i=0}^{[x] }  \sum_{j=1}^{[x] }  j    \big ( \theta(j) \wedge  \theta(i)  \big ) + \sum_{i=0}^{[x] }  \sum_{j \geq [x] }  j    \big ( \theta(j) \wedge  \theta (i)  \big )   \\
 \ll  x  \sum_{i \geq [x]} i  \theta(i)  + \sum_{i=0}^{[x] }  (i+1)^2  \theta(i)  +  \sum_{j \geq 1} j ( j \wedge x ) \theta(j)  \, .
\end{multline*}
Hence
\beq  \label{B5dec6-F}
\sum_{k=1}^n \sum_{i=0}^{{k}} |    \BBE ( {\tilde  \Delta}_{2,3,i}^{(1)}  )  | \ll n  \Big ( 1 + \sum_{k \geq 1} k ( k \wedge x)  \theta(k)  \Big )  \, .
\eeq
With similar arguments, we infer that 
\beq \label{B6dec6-F} 
\sum_{k=1}^n  \big \{  \big |   \BBE \big (  \varphi_{n-k}^{(3)} (\tilde S_{k-1} ) Z_k^3  \big )    \big |   +   \big |   \BBE  (\Delta_{n,k,3}^{(1)} ) \big |  \big \}  \ll n \Theta_2 \, .
\eeq
Starting from \eqref{linddec6-F} and taking into account \eqref{linddec5-F-rests}, \eqref{B5dec6-F} and \eqref{B6dec6-F}, it follows that,  if one can prove that 
\beq \label{laststep-F}
 \sum_{k =1}^n \sum_{i=1}^{k}  \sum_{j=1}^i  \big |  \gamma_i \BBE  \big \{   \varphi^{(3)}_{n-k}  (\tilde S_{k-j -1})   Z_{k-j}  \big \}   \big | \ll n \Theta_1 \Theta_2 ,
\eeq
then 
$$
\sum_{k =1}^n \BBE (   \Delta^*_{n,k} )    \ll n \big ( \Theta_1 \Theta_2 + \sum_{k=1}^n k ( k \wedge x)  \theta (k) \big ) \, .
$$
This last upper bound together with \eqref{sumdelta-F},  \eqref{deltaetoileF}, \eqref{reste1-F} and \eqref{linddec1-F} will end the proof of 
\eqref{aimFuk}. 
The rest of the proof is devoted to the proof of \eqref{laststep-F}. 
Since $ \BBE  \big \{   \varphi^{(3)}_{n-k}  (0)   Z_{k-j}  \big \}  =0$, we have
\begin{multline*} 
\big | \BBE  \big \{   \varphi^{(3)}_{n-k}  (\tilde S_{k-j -1})   Z_{k-j}  \big \}   \big |   
=  \Big |    \sum_{u=1}^{k-j-1}  \BBE  \big \{  \big (   \varphi^{(3)}_{n-k}  (\tilde S_{k-j -u})  -  \varphi^{(3)}_{n-k}  (\tilde S_{k-j -u- 1})   \big )  Z_{k-j}  \big \}  \Big | \\
=   \Big |    \sum_{u=1}^{k-j-1}  \int_0^1  \BBE  \big \{  \varphi^{(4)}_{n-k}    (\tilde S_{k-j -u- 1} +t Z_{k-j-u })  Z_{k-j -u} Z_{k-j}  \big \}   dt  \Big |  
\ll  \sum_{u=1}^{k-j-1} \theta (u)    \, .
\end{multline*}
Therefore
\[
 \sum_{k =1}^n \sum_{i=1}^{k}  \sum_{j=1}^i  \big | \gamma_i  \BBE  \big \{   \varphi^{(3)}_{n-k}  (\tilde S_{k-j -1})   Z_{k-j}  \big \}  \big |    
\ll   \sum_{k =1}^n \sum_{i=1}^{k}  \sum_{j=1}^i   \sum_{u=1}^{k-j-1} \theta (u)   \theta(i)  \
\ll n   \Theta_1 \Theta_2   \, ,
\]
which ends the proof of \eqref{laststep-F}  and then of \eqref{aimFuk}.  This ends the proof of the theorem.  \qed

\section{Proof of Theorems \ref{thSAPfaible} and  \ref{thSAP} (case $\sigma^2 >0$)} \label{SectASIP}

\setcounter{equation}{0}

Starting from the construction of the $V_{k,L}$ given in \eqref{defVkN}, we now construct a suitable sequence $(Z_i)_{i \geq 1}$ of i.i.d.  Gaussian random variables, centered and with variance $\sigma^2$.   Let $Z_1 = \sigma \Phi^{-1} (\delta_1)$. For any
$L \in {\mathbb N}$ and any $k \in \{1, \cdots, 2^{L -m(L)} \}$ 
the random variables $(Z_{2^L + (k-1)2^{m(L)} +1}, \ldots , Z_{2^L + k 2^{m(L)}})$ 
are defined in the following way. If $m(L)=0$, then 
$Z_{2^L + k 2^{m(L)}} = V_{k,L}$.
If $m(L)>0$, then by the Skorohod lemma \cite{Sk}, there exists a measurable function $g$ from ${\mathbb R} \times [0,1]$ 
in ${\mathbb R}^{2^{m(L)}}$ such that, for any pair 
$(V,\delta)$ of independent random variables with respective laws $N(0,\sigma^2 2^{m(L)})$ and 
the uniform distribution over $[0,1]$, $g(V,\delta) = (N_1, \ldots , N_{2^{m(L)}})$ is a Gaussian random vector with i.i.d. components such that $V = N_1 + \cdots + N_{2^{m(L)}}$. 
We then set 
$$
(Z_{2^L + (k-1)2^{m(L)} +1}, \ldots , Z_{2^L + k 2^{m(L)}})  = g (V_{k,L} , \delta_{2^L + (k-1)2^{m(L)} +1} ) \, . 
$$
The so defined sequence $(Z_i)$ has the prescribed distribution. 
\par\medskip
Set $S_j = \sum_{i=1}^j X_i$ and $T_j = \sum_{i=1}^j Z_i$. Let 
$D_L = \sup_{ \ell \leq 2^{L}} | \sum_{i = 2^L +1}^{2^L + \ell} (X_i -Z_i)| $.
Then, proceeding exactly as in \cite{MR12}, page 394,  for any $N$ in $\BBN^*$, 
\begin{eqnarray} \label{1dec} 
\sup_{1 \leq k \leq 2^{N+1} }|S_k -B_k|  \leq 
|X_1 - Z_1| + D_0 + D_1 + \cdots + D_N \, . 
\end{eqnarray}
\par
It remains to bound up the random variables $D_L$. We first notice that the following decomposition is valid: 
\beq \label{decsup} 
D_L  \leq D_{L,1} + D_{L,2} \, ,
 \eeq 
where, recalling that $ U_{k,L} = \sum_{i \in I_{k,L}} X_i$, we set
$$ 
D_{L,1}:= \sup_{k \leq 2^{L-m(L)} } \Big|
\sum_{ \ell =1}^k (U_{\ell,L} - V_{\ell,L}) \Big|
\ \text{and}\ 
D_{L,2}:= \sup_{ k \leq 2^{L-m(L)}} \sup_{ \ell \in I_{k,L} } 
\Big| \sum_{i = \inf I_{k,L}  }^{\ell} (X_i -Z_i) \Big| \, .
$$ 
In order to bound up $D_{L,1}$ and  $D_{L,2} $ we shall use the two lemmas below.

\begin{Lemma} \label{LemDL1} Assume that $\sum_{k \geq 1} k \theta (k) < \infty$. Then there exists a positive constant $C$ such that, for any integer $m(L) $ in $[0,L]$,
\begin{equation} \label{DL1borneL2}
\|D_{L,1}\|^2_2 \leq C 2^{L-m(L) }  \big ( 1 + \sum_{k \geq 1} k ( k \wedge 2^{m(L)/2})  \theta (k)  \big ) \, .
\end{equation}
\end{Lemma}

\begin{Lemma} \label{LemDL2} Assume  that $\sum_{k \geq 1} k \theta (k) < \infty$. Then there exist  positive constants $c$ and $C$ such that, 
for any positive $\lambda$ and any integer $m(L) $ in $[0,L]$,  
\begin{equation}
\BBP  ( D_{L,2} \geq 2\lambda )  \leq  C 2^L \exp \bigl( - 2^{-m(L)} \lambda^2 /c \bigr) + 
C 2^L\lambda^{-4}  \big ( 1 + \sum_{k \geq 1} k ( k \wedge \lambda)  \theta (k)  \big )   \, . \label{apppropineg}
\end{equation}
\end{Lemma}

\noindent
{\bf Proof of Lemma \ref{LemDL1}.} For any $\ell \in \{ 1, \cdots, 2^{L-m(L)} \} $, let
$\widetilde U_{\ell,L} =U_{\ell,L} - \BBE_{2^L + (\ell-1)2^{m(L)}} (U_{\ell,L})$. 
Then $(\widetilde U_{\ell,L})_{\ell \geq 1}$ is a strictly
stationary sequence of martingale differences adapted to the
filtration $({\mathcal F}_{2^L + \ell 2^{m(L)}})_{\ell \geq 1}$. Notice first that
\beq \label{dec0item2}
\|D_{L,1}\|_2 \leq  \Big \Vert \sup_{k \leq 2^{L-m(L)}} \Big  | \sum_{ \ell =1}^k (\widetilde
U_{\ell,L} - V_{\ell,L}) \Big   | \Big  \Vert_2 + \Big \Vert \sup_{k \leq 2^{L-m(L)}}  \Big  | \sum_{ \ell =1}^k (\widetilde
U_{\ell,L} - U_{\ell,L})  \Big   | \Big  \Vert_2  \, .
\eeq
Let us deal with the first term on right hand. Proceeding as in the proof of  Lemma 4.12 in \cite{MR12}, we have 
\[
 \Big \Vert \sup_{k \leq 2^{L-m(L)}}  \Big  | \sum_{ \ell =1}^k (\widetilde
U_{\ell,L} - V_{\ell,L})  \Big   | \Big  \Vert^2_2   \leq 16  \sum_{\ell=1}^{2^{L-m(L)}}
\Vert U_{\ell, L} - V_{\ell, L}\Vert^2_2 \, .
\]
%Since $V_{\ell,L}$ is independent of 
%${\mathcal F}_{2^L + (\ell - 1) 2^{m(L)} }$, the sequence $(\widetilde
%U_{\ell,L} - V_{\ell,L})_\ell$ is a martingale difference sequence with respect to the
%nondecreasing filtration $({\mathcal F}_{2^L + \ell 2^{m(L)}})_\ell$. Hence, by the 
%Doob-Kolmogorov maximal inequality, we get that 
%\begin{eqnarray*}
% \Big \Vert \sup_{k \leq 2^{L-m(L)}}  \Big  | \sum_{ \ell =1}^k (\widetilde
%U_{\ell,L} - V_{\ell,L})  \Big   | \Big  \Vert^2_2 & \leq & 4
%\sum_{\ell=1}^{2^{L-m(L)}} \Vert \tilde  U_{\ell, L} -  V_{\ell, L}
%\Vert^2_2 \\
%& \leq & 8 \sum_{\ell=1}^{2^{L-m(L)}}  \Vert \widetilde
%U_{\ell,L} -  U_{\ell,L} \Vert^2_2 +8 \sum_{\ell=1}^{2^{L-m(L)}}
%\Vert U_{\ell, L} - V_{\ell, L}\Vert^2_2 \, .
%\end{eqnarray*}
%Since $V_{\ell,N}$ is independent of  $\mathcal{F}_{2^L +
%(\ell-1)2^{m(L)} }$, $\BBE_{2^L + (\ell-1)2^{m(L)}} (V_{\ell,L}) =0$.
%Consequently,
%$$
% \Vert \tilde  U_{\ell, L} -  U_{\ell, L} \Vert_2 = \Vert \BBE_{2^L + (\ell-1)2^{m(L)}} (U_{\ell, L} - V_{\ell,L}) \Vert_2
%\leq \Vert U_{\ell, L} - V_{\ell,L}\Vert_2 \, . 
%$$ 
By using Theorem 2.1(b) in \cite{DMR22}, it follows that
\beq \label{maj5DR}
 \Big \Vert \sup_{k \leq 2^{L-m(L)}}  \Big  | \sum_{ \ell =1}^k (\widetilde
U_{\ell,L} - V_{\ell,L})  \Big   | \Big  \Vert^2_2  \ll 2^{L-m(L)}   \Bigl ( 1 + \sum_{k \geq 1} k ( k \wedge 2^{m(L)/2})  \theta (k)  \Bigr ) \, .
\eeq
We  deal now  with the second term in the right hand side of (\ref{dec0item2}). Using Proposition 1 in \cite{DR00}, we obtain that
\begin{eqnarray} \label{inegadedrio}
& &  \Big \Vert \sup_{k \leq 2^{L-m(L)}}  \Big  | \sum_{ \ell =1}^k (\widetilde
U_{\ell,L} - U_{\ell,L})  \Big   | \Big  \Vert^2_2   \leq  4 \sum_{k=1}^{2^{L-m(L)}} \Vert\BBE_{2^L + (k-1)2^{m(L)}} (U_{k,L}) \Vert_2^2
\nonumber \\
& & \quad \quad + 8 \sum_{k=1}^{2^{L-m(L)} -1 }
  \Big \Vert\BBE_{2^L + (k-1)2^{m(L)}} (U_{k,L}) \big ( \sum_{i= k +1}^{2^{L-m(L)}}\BBE_{2^L + (k-1)2^{m(L)}} (U_{i,L}) \big ) \Big  \Vert_1 \, .
\end{eqnarray}
Stationarity leads to
\[
\Vert\BBE_{2^L + (k-1)2^{m(L)}} (U_{k,L}) \Vert_2^2= \Vert\BBE_{0} (S_{2^{m(L)}}) \Vert_2^2  
\leq 2 \sum_{i=1}^{2^{m(L)}} \sum_{j=1}^i {\mathbb E} | X_j {\mathbb E}_0 (X_i)|  \leq 2 \sum_{i=1}^{2^{m(L)}} i  \theta (i)  \ll 1 \, .
\]
Consequently, 
 \beq \label{maj2DR}
\sum_{k=1}^{2^{L-m(L)}} \Vert\BBE_{2^L + (k-1)2^{m(L)}} (U_{k,L}) \Vert_2^2  \ll  2^{L-m(L)}  \, .
\eeq
We now bound up the second term in the right hand side of (\ref{inegadedrio}). Stationarity yields 
$$
 \Big  \Vert \BBE_{2^L + (k-1)2^{m(L)}} (U_{k,L}) \Big ( \sum_{i= k +1}^{2^{L-m(L)}}\BBE_{2^L + (k-1)2^{m(L)}} (U_{i,L}) \Big )  \Big \Vert_1  \leq  \sum_{j=1}^{2^{m(L)}} 
\sum_{i=2^{m(L)}+ 1}^{2^L - (k-1) 2^{m(L)}} {\mathbb E} | X_j {\mathbb E}_0 (X_i)| \, .
$$
Therefore
\begin{equation} \label{maj3DR}
\sum_{k=1}^{2^{L-m(L)} -1 }
 \Big \Vert \BBE_{2^L + (k-1)2^{m(L)}} (U_{k,L}) \Big ( \sum_{i= k +1}^{2^{L-m(L)}}\BBE_{2^L + (k-1)2^{m(L)}} (U_{i,L}) \Big )  \Big \Vert_1  
  \ll 2^{L-m(L)} 
\sum_{i=2^{m(L)}+ 1}^{2^L }  i  \theta (i)  \, .
\end{equation}
Starting from (\ref{inegadedrio}) and considering the bounds (\ref{maj2DR}) and (\ref{maj3DR}), 
we get that
\begin{equation} \label{maj4DR}    \Big   \Vert \sup_{k \leq 2^{L-m(L)}}  \Big  | \sum_{ \ell =1}^k (\widetilde
U_{\ell,L} - U_{\ell,L})   \Big  |  \Big \Vert_2^2  \ll 2^{L-m(L)} \, .
\end{equation}
Starting from (\ref{dec0item2}) and  considering the bounds (\ref{maj5DR}) and 
(\ref{maj4DR}), we then get (\ref{DL1borneL2}),
which ends the proof of Lemma \ref{LemDL1}. \qed
\par\medskip

\noindent
{\bf Proof of Lemma \ref{LemDL2}.}  It follows the lines of the proof of \cite[Lemma 4.1]{MR12} with the difference that  Theorem 
\ref{Fuk-inequality-theta} is used instead of \cite[Proposition A.2]{MR12}. 
%By the triangle inequality together with the stationarity of the sequences $(X_i)_i$ and $(Z_i)_i$,  for any positive $\lambda$,
%\beq \label{decdl2} {\mathbb P}
%( D_{L,2} \geq 2 \lambda ) \leq 2^{L-m(L)} {\mathbb P} \Bigl(  \sup_{
%\ell \leq 2^{m(L)}} | S_{\ell} | \geq \lambda \Bigr) +  2^{L-m(L)} {\mathbb P} 
%\Bigl(  \sup_{\ell \leq 2^{m(L)}} | T_{\ell} | \geq \lambda
%\Bigr)  \, .\eeq
%By L\'evy's inequality (see for instance Proposition 2.3 in Ledoux and Talagrand (1991)),
%\begin{equation} {\mathbb P} \Bigl(  \sup_{
% \ell\leq 2^{m(L)}} | T_{\ell} | \geq \lambda
%\Bigr) \leq  2 \exp  \Bigl( - \frac{\lambda^2}{2 \sigma^2 2^{ m(L)}
%} \Bigr)   \, . \label{LI}
%\end{equation}
%On the other hand, applying Theorem \ref{Fuk-inequality-theta} , we get that
%\[
%\BBP  \Bigl(  \sup_{ \ell \leq 2^{m(L)}} | S_{\ell} | \geq \lambda \Bigr)  \leq  
%c_1\exp \Bigl( - \frac{\lambda^2}{c_2
%2^{m(L)}} \Bigr) + c_3 2^{m(L)}\lambda^{-4} \big ( 1 + \sum_{k \geq 1} k ( k \wedge \lambda)  \theta (k)  \big )    \, .
%\]
%Collecting the above inequalities, we then get Lemma \ref{LemDL2}. 
\qed
\par\medskip

\noindent{\bf End of the proof of Theorem \ref{thSAPfaible}.}  Let us start by proving Item (i). In case $p \in ]2,3[ $ with the condition $ \theta (k) = O ( k^{1-p} )$  and in case $p=3$ with the condition 
$\sum_{k >0} k \theta (k) < \infty$ the result can be proved exactly as in  \cite{MR12} (indeed, in the bounded case, these authors  could have used the coefficient $\theta (k)$ instead of their coefficient $\alpha_{2, {\mathbf X}} (k)$). 
Now, we turn to the case $p \in ]3,4] $. Note that under the condition $ \theta (k) = O ( k^{1-p} )$ if  $p \in ]3,4[ $ and  $\sum_{k >0} k^2 \theta (k) < \infty$ if $p=4$, 
\beq \label{factobviousseries}
 \sum_{k \geq 1} k ( k \wedge \lambda)  \theta (k)  \leq C  ( 1+   \lambda^{4-p }  )  \, .
\eeq
Therefore, by Lemma \ref{LemDL2}, simple computations lead to 
\begin{equation} \label{DL2carre}
\BBE ( D_{L,2}^2 ) = 2 \int_0^{\infty} x   {\mathbb P}
( D_{L,2} \geq x  )   d x   \leq C  \big (  L 2^{m(L)} +   2^{2L/p}  \big )  \, ,
\end{equation}
and,  by Lemma \ref{LemDL1}, 
\begin{equation} \label{DL1carre}
\BBE ( D_{L,1}^2 )    \leq C  2^L 2^{ ( 1-p/2) m(L)}  \, .
\end{equation}
Choosing 
\beq \label{selectionML}
  m(L) = [ 2(L -   \log_2 L)/p ] \ \text{so that}\ 
2^{ -1 + 2L/p} L^{ -2/p} \leq 2^{m(L)} \leq 2^{2L/p} L^{ -2/p} \, , 
\eeq
for the construction of the Gaussian sequence, 
Item (i) follows (above square brackets designate as usual the integer part and $\log_2 (x)= (\log x)/(\log 2)$).  

We turn now to Item (ii). Let us complete the proof when $p=4$, meaning that  $ \theta (k) = O ( k^{-3} )$. In this case,  
\beq \label{factobviousseriesp=4}
 \sum_{k \geq 1} k ( k \wedge \lambda)  \theta (k)  \leq C  ( 1+   \log \lambda   )  \, , 
\eeq
for any $\lambda \geq 1$, 
and then the term $ 2^{2L/p} $ appearing in \eqref{DL2carre} has to be replaced by 
$(1+m(L) ) 2^{m(L)}$ and the right-hand side of \eqref{DL1carre} will be $C  L 2^{L-m(L)}$. Choosing  $m(L) = [L/2]$ completes the proof.  Finally, the case $p=3$ with the condition $ \theta (k) = O ( k^{-2} )$ can be handled similarly  
by taking into account Lemmas 4.2 and 4.1  in  \cite{MR12}  instead of our Lemmas \ref{LemDL1} and \ref{LemDL2}. 

\medskip

\noindent{\bf End of the proof of Theorem \ref{thSAP}.}   In case $p \in ]2,3[ $ and $ \sum_{k>0} k^{p-2} \theta(k)< \infty$  the result can be proved exactly as in Theorem 2.1 Item 1a) in  \cite{MR12} (using the coefficient $\theta(k)$ instead of their coefficient $\alpha_{2, {\mathbf X}} (k)$).  Similarly, in case $p \in ]2,3[ $ and $ \theta(k) = O ( k^{1-p} )$ or $p=3$ and $\sum_{k >0}k \theta(k) < \infty$,  the result can be proved exactly as in Theorem 2.1 Item 1b) in  \cite{MR12} (with $\theta(k)$ instead of  $\alpha_{2, {\mathbf X}} (k)$). 

Let us now complete  the proof of Item (i) when $p \in ]3,4[$ and $ \sum_{k>0} k^{p-2} \theta(k)< \infty$. In this case we select $m(L)$ as in \eqref{selectionML} and set 
\beq \label{deflambda}
\lambda_L = \kappa 2^{m(L)/2} \sqrt{L} \, , 
\eeq
with $\kappa= \sqrt{2 c   \log 2} $ where $c$ is the positive constant of Lemma \ref{LemDL2}.  For this choice,  
\beq \label{cons1selectionlambda}
\sum_{L > 0} 2^L \exp \bigl( -\lambda_L^2 2^{-m(L)}/c \bigr) 
= \sum_{L \geq 0} 2^{L-2L} < \infty
\ \text{and}\ 
\sum_{L > 0} 2^{L} \lambda_L^{-4} < \infty \, .
\eeq
In addition, using the fact that  $\sum_{k \geq 1} k ( k \wedge (a \lambda) ) \theta(k) \leq a  \sum_{k \geq 1} k ( k \wedge \lambda ) \theta(k)  $, for any $a \geq 1$ and any positive  $\lambda$, we get, since $ \sum_{k>0} k^{p-2} \theta(k)< \infty$, 
\[
\sum_{L > 0}  2^L \lambda_L^{-4}   \sum_{k \geq 1} k ( k \wedge  \lambda_L ) \theta(k)  
 \leq  \sum_{L > 0}  2^L L^{1/2} \lambda_L^{-4}   \sum_{k \geq 1} k ( k \wedge  2^{m(L)/2} ) \theta(k)  
 < \infty \, .
\]
Therefore Lemma \ref{LemDL2}  entails that $\sum_{L > 0} {\mathbb P} (
D_{L,2} \geq 2 \lambda_L ) < \infty $ implying, via the Borel-Cantelli lemma, that 
\begin{equation}
D_{L,2} = O ( 2^{L/p} L^{1/2-1/p}) 
\text{  a.s. } \label{DL2}
\end{equation}
On another hand,  from (\ref{DL1borneL2}) together with the Markov inequality, 
$$
\sum_{L>0} {\mathbb P} ( D_{L,1} \geq \lambda_L) \leq 
C \sum_{L>0}  2^L L^{-1} 2^{-2m(L)}   \sum_{k \geq 1} k ( k \wedge  2^{m(L)/2} ) \theta(k)   < \infty\, ,
$$ 
since $ \sum_{k>0} k^{p-2} \theta(k)< \infty$. Hence, by the Borel-Cantelli
lemma, 
\beq
D_{L,1} = O ( 2^{L/p} L^{1/2-1/p}) 
\text{  a.s. } \label{DL1}
\eeq
Finally Item (i) when $p \in ]3,4[$  follows from  (\ref{1dec}), (\ref{decsup}),  (\ref{DL2}) and (\ref{DL1}).

\smallskip

We complete now  the proof of Item (ii) when $p \in ]3,4[$ and $ \theta(k) = O (k^{1-p})$ or when $p=4$ and $\sum_{k>0}k^2 \theta(k) < \infty$. In these cases we select $m(L)$ as follows: let $\varepsilon >0$ and set 
\beq \label{selectionMLbis}
  m(L) = [ 2(L + \varepsilon \log_2 L )/p] \ \text{so that}\ 
  2^{ -1+2L/p} L^{ 2 \varepsilon /p} \leq 2^{m(L)} \leq 2^{2L/p} L^{ 2 \varepsilon /p} \, . 
\eeq
We still define $\lambda_L$ by \eqref{deflambda}.  
For this choice of $\lambda_L$, the convergences in \eqref{cons1selectionlambda} still hold. 
In addition, taking into account \eqref{factobviousseries}, for any $p \in ]3,4]$, under the conditions on 
$(\theta(k))_{k>0}$, we get 
\beq \label{BorneD1}
\sum_{L > 0}  2^{L} \lambda_L^{-4}   \sum_{k \geq 1} k ( k \wedge  \lambda_L ) \theta(k)  
\leq C  \sum_{L > 0}  2^L  \lambda_L^{-p}   
 < \infty \, .
\eeq
Therefore Lemma \ref{LemDL2}  entails that $\sum_{L > 0} {\mathbb P} (
D_{L,2} \geq 2 \lambda_L ) < \infty $ implying, via the Borel-Cantelli lemma, that 
\begin{equation}
D_{L,2} = O ( \lambda_L) =  O ( 2^{L/p} L^{1/2+ \varepsilon/p}) 
\text{  a.s. } \label{DL2-2}
\end{equation}
On another hand,  from (\ref{DL1borneL2}) and \eqref{factobviousseries} together with the Markov inequality, 
\beq \label{BorneD2}
\sum_{L>0} {\mathbb P} ( D_{L,1} \geq \lambda_L) \leq 
C \sum_{L>0}  \bigl( 2^L L^{-1} 2^{-2m(L)} \bigr) 2^{m(L) ( 2- p/2)}    < \infty\, .
\eeq
Hence, by the Borel-Cantelli
lemma, 
\beq
D_{L,1} = O ( \lambda_L) =  O ( 2^{L/p} L^{1/2+ \varepsilon/p}) 
\text{  a.s. } \label{DL1-2}
\eeq
Finally Item (ii) when $p \in ]3,4]$  follows from  (\ref{1dec}), (\ref{decsup}),  (\ref{DL2-2}) and (\ref{DL1-2}). 

\smallskip

We turn now to the proof of Item (iii).  When $p=4$, meaning that  $\theta (k)  = O( k^{-3} )$, instead of  \eqref{factobviousseries}, we use \eqref{factobviousseriesp=4}. We still select $\lambda_L$ by \eqref{deflambda} but with the following choice of $m(L)$: for $\varepsilon >0$
\beq \label{choixdemLp}
  m(L) = [ 2 (L+  (1+ \varepsilon ) \log_2 L )/p ] \ \text{so that}\ 
 2^{-1+ 2L/p}  L^{ 2(1+\varepsilon) /p} \leq 2^{m(L)} \leq 2^{ 2 L/p} L^{2 (1+\varepsilon) /p} \, . 
\eeq
The computations \eqref{BorneD1} and \eqref{BorneD2}  are then replaced by the following ones: 
\[
\sum_{L > 0}  \frac{2^{L} }{\lambda_L^{4} } \sum_{k>0}  ( k^2 \wedge  k \lambda_L ) \theta(k)  
\leq C  \sum_{L > 0}  \frac{2^{L}  }{\lambda_L^{4} }   \log \lambda_L 
 \, ,\ 
%\]
%and
%\[
\sum_{L>0} {\mathbb P} ( D_{L,1} \geq \lambda_L) \leq 
C \sum_{L>0}  \frac{  L 2^L }{2^{m(L)} \lambda_L^2 }   \, .
\]
Note that the above upper bounds are finite for these selections of $m(L)$ and $\lambda_L$. 
The rest of the proof is unchanged compared to the previous cases. 

\smallskip

It remains to prove Item (iii) when $p=3$, meaning that  $\theta (k)  = O( k^{-2})$. The differences with the case $p=4$ are that  \cite[Lemma 4.1]{MR12} is used instead of our Lemma  \ref{LemDL1} and  $ \Vert D_{L,1} \Vert_2^2 \leq C 2^{L-m(L)/2} m(L)$. This upper bound on $\Vert D_{L,1} \Vert_2^2$ comes from a slight modification of the proof of \cite[Lemma 4.2]{MR12}, taking into account that
 $\sum_{k \geq 1}  ( k \wedge \lambda)  \theta (k)  \leq C  ( 1+  \log \lambda  )$ for any 
 $\lambda \geq 1$. 
In addition,  $m(L)$ is selected by \eqref{choixdemLp} with $p=3$. \qed

\section{The degenerate case} \label{Sectiondegenerate}

\setcounter{equation}{0}

In all this section, we shall denote $ \theta_{X,1,1}(k) $ by $ \theta(k)$ for all $k \in {\mathbb N}$.

\begin{Proposition} \label{Thdegenerate} Assume that $\Vert X_0 \Vert_{\infty} \leq M < \infty$ and $ \sum_{k \geq 0} \theta(k)  < \infty$. Suppose  in addition  that  $\sigma^2=0  $. 
Then, for any $q \geq 1$, 
\[
\BBE ( |S_n|^{q} )  \leq q  (2M )^{q} \sum_{k \geq 0} ( k+1)^{q-1} \theta(k) \, .
\]
%\[
%\BBE ( |S_n|^p ) \leq 4^{p-1} n \int_{0}^1 R^{p-1} (u) Q(u) du 
%\]
%and
%\[
%\BBE ( |S_n|^{p-1} ) \leq 4^{p-2}  \int_{0}^1 R^{p-1} (u)  du \, .
%\]
\end{Proposition}

\begin{Theorem}   \label{CordegenerateL2}    Let  $p > 2$ and $S_n^* = \max_{1 \leq k \leq n} |S_k|$.   Assume that $\Vert X_0 \Vert_{\infty} \ < \infty$ and  $ \theta (k) = O ( k^{1-p} )$. Suppose in addition that   $\sigma^2=0  $.  Then, for any $r \in [1, p[ $, $ \Vert S_n^* \Vert_r = O (n^{1/p})$.
\end{Theorem}
\begin{Remark}
Since for any increasing sequence $(b_n)_{n \geq 1}$, $\sum_{n >0} n^{-1} \BBP (  S_n^* > b_n ) < \infty$ implies that $S_n = o ( b_n) $ almost surely, it follows that,  under the assumptions of the theorem, $S_n = o (n^{1/p}  (\log n)^{  \varepsilon + 1/p}) $ almost surely for any $\varepsilon >0$. 
\end{Remark}

\begin{Theorem}   \label{Cordegenerate}  Let  $p \geq 2$ and $S_n^* = \max_{1 \leq k \leq n} |S_k|$.  Assume that $\Vert X_0 \Vert_{\infty} \ < \infty$ and  $\ \sum_{k \geq 0} ( k+1)^{p-2} \theta(k) < \infty$. Suppose in addition that   $\sigma^2=0  $.
Then, for any $\alpha \in ]0, 1 [$ and any  $\varepsilon >0$, $\sum_{n >0} n^{\alpha p-2} {\mathbb P} ( S_n^* > \varepsilon n^{\alpha}  ) < \infty$. 
% \beq \label{BK}
%\sum_{n >0} n^{\alpha p-2} {\mathbb P} ( S_n^* > \varepsilon n^{\alpha}  ) < \infty \, .
%\eeq
Consequently $S_n = o (n^{1/p} ) $ almost surely. 
\end{Theorem}

\noindent \textbf{Proof of Proposition \ref{Thdegenerate} }.  We start the proof with the following lemma. 
\begin{Lemma}  \label{Lemmadec} Assume that $X_0 \in {\mathbb L}^2$, 
\begin{itemize}
\item[(a)] $\BBE_0 (S_n)$ converges in ${\mathbb L}^1$,
\item[(b)]  $ \lim_{n \rightarrow \infty}\BBE (X_0 (X_0 + 2 S_n)) =0$,
\item[(c)]  $ \limsup_{N \rightarrow \infty} \limsup_{n \rightarrow \infty}\BBE ( \BBE_{-N} (X_0 ) (X_0 + 2 S_n)) \leq 0$. 
\end{itemize}
Then, for any integer $i$,  $X_i = g_{i-1} - g_i $ almost surely where $g_i = \sum_{k \geq i  +1} \BBE_{i} ( X_k) $. 
\end{Lemma}
\noindent {\bf Proof of Lemma \ref{Lemmadec}.}   Let $N$ be a fixed positive integer. Set $d_{i,N} = \sum_{k=0}^{N-1}  P_i (X_{k+i}) $ where $P_i =  \BBE_i   - \BBE_{i-1} $. 
Define also $g_{i,N} =  \sum_{k =1}^N \BBE_{i} ( X_{k+i})  $ and $Y_{i,N} =  \BBE_{i} ( X_{N+i})  $.  With these notations, the following decomposition is valid:
\begin{equation} \label{deccobwithN}
X_i = d_{i,N }  + g_{i-1,N} - g_{i,N}   + Y_{i,N} \, .
\end{equation}
Since $N$ is fixed, all the random variables in the above decomposition are in ${\mathbb L}^2$. Moreover, from the fact that $(d_{i,N})_{i \in {\mathbb Z}}$  is a stationary sequence of martingale differences, 
\begin{equation} \label{dconsitemb0}
\Vert d_{0,N } \Vert_1 \leq  \Vert d_{0,N } \Vert_2  \leq \limsup_{n \rightarrow \infty } n^{-1/2} \big \Vert  \sum_{i=1}^n d_{i,N} \big \Vert_2 \, .
\end{equation}
Next, from \eqref{deccobwithN},
$\sum_{i=1}^n d_{i,N}  = S_n + g_{n,N} - g_{0,N} - \sum_{i=1}^n Y_{i,N}$,
which implies that 
\[
 n^{-1/2} \big \Vert  \sum_{i=1}^n d_{i,N} \big \Vert_2 \leq   n^{-1/2}  \Vert  S_n  \Vert_2 +  2  n^{-1/2}  \Vert  g_{0,N} \Vert_2 +   n^{-1/2} \big \Vert  \sum_{i=1}^n Y_{i,N} \big \Vert_2 \, .
\]
Now, by item (b), $ \lim_{n \rightarrow \infty}  n^{-1/2}  \Vert  S_n  \Vert_2 =0$. Hence
\begin{equation} \label{dconsitemb1}
\limsup_{n \rightarrow \infty} n^{-1/2} \big \Vert  \sum_{i=1}^n d_{i,N} \big \Vert_2 \leq  \limsup_{n \rightarrow \infty}  n^{-1/2} \big \Vert  \sum_{i=1}^n Y_{i,N} \big \Vert_2 \, .
\end{equation}
Next, by stationarity and the properties of the conditional expectation, we infer that  
\[
 \big \Vert  \sum_{i=1}^n Y_{i,N} \big \Vert^2_2  =  \sum_{k=0}^{n-1} \BBE \big ( \BBE_{-N} (X_0) ( X_0  + 2  S_{k}   ) \big ) \, .
\]
Now Item (c) combined with the Cesaro Lemma entails that 
\begin{equation} \label{dconsitemb2}
 \lim_{N \rightarrow \infty}\limsup_{n \rightarrow \infty} n^{-1}   \big \Vert  \sum_{i=1}^n Y_{i,N} \big \Vert^2_2  = 0 \, .
\end{equation}
Taking into account \eqref{dconsitemb0}-\eqref{dconsitemb2}, we derive that $d_{i,N}$ converges to $0$ in ${\mathbb L}^1$ as $N $ tends to $ \infty$.  Therefore the lemma follows by taking into account the decomposition \eqref{deccobwithN} and noting that, by item (a), 
$g_{i,N}$ converges to $g_i$ in ${\mathbb L}^1$ and $Y_{i,N}$ converges to $0$ in ${\mathbb L}^1$. \qed

\medskip   

We now prove Proposition \ref{Thdegenerate}.   We start by noticing that $\sum_{k \geq 1} \theta ( k) < \infty$ and  
$\sigma^2=0$  imply items (a), (b) and (c) of Lemma  \ref{Lemmadec}. More precisely, to prove item (c), starting from Inequality (4.7) in \cite{Rio17}, we obtain that 
$\big |  \BBE ( \BBE_{-N} (X_0 )X_k )  \big |  \leq M \theta( k+N)$,
implying that 
$\big |  \BBE ( \BBE_{-N} (X_0 ) (X_0 + 2 S_n)  )  \big |  \leq 2M  \sum_{\ell \geq N} \theta (\ell )$.
So, overall, under the conditions of Theorem \ref{Thdegenerate}, 
\beq \label{deccobord}
S_n = g_0 - g_n  \mbox{ in  ${\mathbb L}^1$ where } g_m = \sum_{k >0} \BBE ( X_{k+m} | {\mathcal F}_m )  \, .
\eeq

Now, let $q \geq 1$ and $r $ in $ ]1, + \infty]$ be the conjugate exponent, that is $1/q + 1/r =1$. Note first that,
by the Riesz-Fisher theorem,  
$\Vert g_0 \Vert_q =  \sup  \big \{  \BBE ( g_0 Z)   \, : \, Z \in {\mathcal L}_{r}  \big \}$, 
where  $ {\mathcal L}_{r}$ is the class of nonnegative random variables $Z$ such that $   \Vert Z \Vert_r =1$. But, by \cite[Proposition 1]{DD}, 
\[
 \big |  \BBE ( g_0 Z)  \big | \leq  \sum_{k >0}  |  \BBE \big (  Z X_k   \big )  | \leq M  \sum_{k >0} \int_0^{\theta (k) } Q_Z (u) du  \leq M \int_0^1  \theta^{-1} (u)  Q_Z (u) du  \, ,
\]
where $Q_Z$ is the generalized inverse of $t \mapsto {\mathbb P} (Z >t)$ and $ \theta^{-1} (u) = \sum_{i \geq 0} {\bf 1}_{u < \theta (i)}$. 
Therefore, by H\"older's inequality, 
\[
 \big |  \BBE ( g_0 Z)  \big | \leq M  \Big ( \int_0^1  \big (  \theta^{-1} (u)  \big )^{q} du \Big )^{1/q} \, .
\]
Hence, using the coboundary decomposition \eqref{deccobord} and the stationarity of $(g_m)_{m \in {\mathbb Z}}$, 
\[
\Vert S_n \Vert_{q} \leq 2 \Vert g_0 \Vert_{q} \leq 2 M \Big ( \int_0^1  \big (  \theta^{-1} (u)  \big )^{q} du \Big )^{1/q}     \, .
\]
Now, by inequality (C.5), page 184, in \cite{Rio17}, $\int_0^1  \big (  \theta^{-1} (u)  \big )^{q} du \leq q \sum_{k \geq 0} (k+1)^{q-1} \theta(k)$.   \qed

\smallskip

\noindent \textbf{Proof of Theorem  \ref{CordegenerateL2}}.  Let $q = \max ( 1, r(p-1)/p)$.  We  apply inequality \eqref{1F} with $\varphi(x) = x^q$. Since 
$ \sum_{k>0}  k^{q-1}  \theta (k) < \infty $, Proposition  \ref{Thdegenerate} ensures that $ \Vert S_n \Vert_q \ll 1$.   Therefore, for any $x>0$, 
\beq \label{FNAppli}
\BBP \big ( S_n^* \geq  x  \big ) \ll  x^{-q}  + \min ( 1 , n x^{-p}) \, .
\eeq 
Next, since $ S_n^*  \leq nM$, 
$\Vert S_n^* \Vert_r^r = r \int_0^{nM} x^{r-1} \BBP \big ( S_n^* \geq  x  \big ) dx$.
Hence,  applying \eqref{FNAppli} and taking into account the selection of $q$, we get 
\[
 \Vert S_n^* \Vert_r^r \ll \int_0^{nM} x^{r-1-q} dx +   \int_0^{n^{1/p}} x^{r-1} dx + n   \int_{n^{1/p}}^{\infty} x^{r-1-p} dx \ll n^{r/p} \, . \qed
\]

\smallskip

\noindent \textbf{Proof of Theorem \ref{Cordegenerate} }.  
%The fact that $S_n = o (n^{1/p} ) $ almost surely follows from \eqref{BK} applied with $\alpha = 1/p$.  We now prove \eqref{BK}. 
We start  from inequality \eqref{1F} applied with $\varphi(x) = x^{p-1}$. By Proposition  \ref{Thdegenerate},  $\BBE ( \varphi (S_n) ) \ll 1$. Therefore
\[
\BBP \big ( S_n^* \geq 4 \varepsilon n^{\alpha} \big ) \ll   (\varepsilon n^{\alpha})^{1-p}+ n^{1-\alpha} \varepsilon^{-1} \theta  (  [ \varepsilon n^{\alpha}  ])   \, .
\]
Hence
\[
\sum_{n >0} n^{\alpha p -2}  \BBP \big ( S_n^* \geq 4 \varepsilon n^{\alpha} \big )  \ll  \varepsilon^{1-p} \sum_{n >0} n^{\alpha - 2}  +  \varepsilon^{-1}
    \sum_{n  \geq 1} n^{\alpha  (p -1) - 1}   \theta  (  [ \varepsilon n^{\alpha}  ])   \, .
\]
The first series converges if $\alpha < 1$ and the second one also converges as soon as $\ \sum_{k > 0}k^{p-2} \theta(k) < \infty$.  This ends the proof of the theorem.  \qed

 %%%%%%%%%%%%%%%%%%%%%%%%%%%%%%% BIBLIOGRAPHY %%%%%%%%%%%%%%%%%%%%%%

\end{document}